\documentclass[12pt]{amsart}
\usepackage{etex} \reserveinserts{28} 
\usepackage{color}
\usepackage{amssymb}
\usepackage[all]{xy}
\usepackage{geometry}
\input prepictex
\input pictexwd
\input postpictex

\newtheorem{theorem}{Theorem}[section]
\newtheorem{lemma}[theorem]{Lemma}
\newtheorem{prop}[theorem]{Proposition}
\newtheorem{corollary}[theorem]{Corollary}
\theoremstyle{definition}
\newtheorem{defn}[theorem]{Definition}
\newtheorem{remark}[theorem]{Remark}
\newtheorem{example}[theorem]{Example}

\newcommand\bbC{\mathbb{C}}

\newcommand\bbN{\mathbb{N}}

\newcommand\bbZ{\mathbb{Z}}

\newcommand\bx{\mathbf{x}}

\newcommand\calL{\mathcal{L}}
\newcommand\calM{\mathcal{M}}
\newcommand\calN{\mathcal{N}}
\newcommand\calO{\mathcal{O}}
\newcommand\calP{\mathcal{P}}
\newcommand\calS{\mathcal{S}}
\newcommand\calT{\mathcal{T}}
\newcommand\calU{\mathcal{U}}

\newcommand\calW{\mathcal{W}}

\newcommand\fS{\mathfrak{S}}

\newcommand\ch{\mathrm{ch}}
\newcommand\Frob{\mathrm{Frob}}
\newcommand\Hom{\mathrm{Hom}}
\newcommand\im{\mathrm{im\,}}
\newcommand\Ind{\mathrm{Ind}}
\newcommand\Res{\mathrm{Res}}
\newcommand\spn{\mathrm{span}}

\newcommand{\largewedge}{\mbox{\Large $\wedge$}}

\title{A categorification of the chromatic symmetric function}
\author{Radmila Sazdanovi\'c}
\address[Radmila Sazdanovi\'c]{North Carolina State University\\
	Raleigh, NC 27695}
\email{rsazdanovic@math.ncsu.edu}
\author{	Martha Yip}
\address[Martha Yip]{University of Kentucky\\
	Lexington, KY 40506}
\email{martha.yip@uky.edu}
\begin{document}
\begin{abstract}
The Stanley chromatic symmetric function $X_G$ of a graph $G$ is a symmetric function generalization of the chromatic polynomial, and has interesting combinatorial properties.  We apply the ideas of Khovanov homology to construct a homology of graded $\fS_n$-modules, whose graded Frobenius series $\Frob_G(q,t)$ reduces to the chromatic symmetric function at $q=t=1$.  This homology can be thought of as a categorification of the chromatic symmetric function, and provides  a homological analogue of several familiar properties of $X_G$. 
In particular, the decomposition formula for $X_G$ discovered recently by Orellana and Scott, and Guay-Paquet is lifted to a long exact sequence in homology. 

\end{abstract}
\maketitle 
\parskip=5pt

\section{Introduction}\label{sec.intro}
Let $G$ be a graph with vertex set $V(G)=\{v_1,\ldots, v_n\}$ and edge set $E(G)$.
A {\em proper colouring} of $G$ is a function $\kappa:V(G)\rightarrow \bbN$ such that $\kappa(v_i)\neq \kappa(v_j)$ if an edge is incident to both $v_i$ and $v_j$.
The {\em chromatic polynomial $\chi_G(x)$ of $G$} is a polynomial such that for every $k\in \bbN$, $\chi_G(k)$ is the number of proper colourings of $G$ with $k$ colours.  

\begin{defn}
The {\em chromatic symmetric function of $G$} is defined to be
\begin{equation}
X_G(\bx) = X_G(x_1,x_2,\ldots) = \sum_{\kappa} x_{\kappa(v_1)}\cdots x_{\kappa(v_n)},
\end{equation}
where the sum is over all proper colourings $\kappa:V(G)\rightarrow \bbN$ of $G$.
\end{defn}
The polynomial $X_G$ is a generalization of the chromatic polynomial $\chi_G$ in the sense that
$X_G(1^k) = \chi_G(k)$~\cite[Proposition 2.2]{Sta95}, for $k\in \bbN$.

For each $i\in \bbN$, $\kappa^{-1}(i)$ is an independent set in $G$, so any permutation of $\bbN$ that fixes all but finitely many elements gives another proper colouring of $G$, and therefore $X_G$ is a symmetric function, homogeneous of degree $n$.
Following the standard notation used by Macdonald, let $p_\lambda$ and 
$s_\lambda$ respectively denote the power sum symmetric functions and the Schur symmetric functions.
Given a subset of edges $F\subseteq E(G)$, its {\em partition type} $\lambda(F)$ is the partition associated to the sizes of the connected components of the subgraph of $G$ induced by the edge set $F$.  The following formula~\cite[Theorem 2.5]{Sta95} can be proved by an inclusion-exclusion argument.
\begin{equation}\label{thm.X_in_p}
X_G = \sum_{F\subseteq E(G)} (-1)^{|F|} p_{\lambda(F)}.
\end{equation}
This formula forms the basis of our categorification process.

Categorification can be thought as a way of replacing an $n$-category by an $(n+1)$-category; for example, lifting the Euler characteristic of a topological space to its homology. One of the most successful recent examples of categorification include link homology~\cite{khovanov2000}, which is a new kind of link invariant that lifts the properties of the Jones polynomial and carries a rich algebraic structure.  In this theory, every link is assigned bigraded homology groups whose Euler characteristic is the Jones polynomial, and additionally, link cobordisms are assigned homomorphisms of homology groups.  This categorification has been successfully used in determining topological properties of knots, and gives a purely combinatorial proof~\cite{ras2010} of the Milnor conjecture, also known as Kronheimer-Mrowka theorem.

Chromatic graph homology, developed by~\cite{HR05}, is one of several categorifications of polynomial graph invariants. The construction follows that of Khovanov; a bigraded homology theory  is associated to a graph and a commutative graded algebra in a way that its graded Euler characteristic is the chromatic polynomial evaluated at the $q$-dimension of the algebra. 
There are other categorifications of the chromatic polynomial, and all of them possess a long exact sequence of homology that lifts the deletion-contraction formula for the chromatic polynomial, in the same way that Khovanov homology lifts the skein relations for the Jones polynomial.

In this paper we apply a Khovanov-type construction to the chromatic symmetric function $X_G$.  This process is described in Section~\ref{sec.construction}.  Every spanning subgraph of $G$ is assigned a graded $\mathfrak{S}_n$-module, leading to a chain complex whose differential maps are defined based on the Boolean lattice structure of the set of spanning subgraphs of $G$.  The Frobenius series $\mathrm{Frob}_G(q,t)$ of the resulting bigraded homology, which we call the chromatic symmetric homology, specializes to the chromatic symmetric function, naturally expressed in the Schur basis. 

In Section~\ref{sec.properties}, we give the analogue of several familiar properties of the chromatic symmetric function on the categorified level.
In particular, we consider the homology of graphs which contain a loop or multiple edges, and the homology of a disjoint union of graphs.

Unlike the chromatic polynomial $\chi_G$, the chromatic symmetric function $X_G$ does not satisfy a deletion-contraction formula, but as shown by~\cite{os2014} and~\cite{guaypaquet2013}, if $G$ contains a $3$-cycle, then $X_G$ satisfies a
recursive formula involving the deletion of two of the edges of the $3$-cycle.  In Section~\ref{sec.mv}, we show how this recursive formula for $X_G$ can be lifted to a long exact sequence in homology.

To conclude, we provide computations of homology for small graphs in Section~\ref{sec.computations}.  Finally, it is our hope that chromatic symmetric homology will be useful in addressing open problems regarding the chromatic symmetric function, such as the characterization of graphs whose $X_G$ is Schur-positive (\cite{gasharov1996}), and whether $X_G$ is a complete invariant for trees (\cite{mmw2008}).  Also, we hope that this construction may shed light on the chromatic quasisymmetric function of~\cite{sw2014}.


\subsection{Preliminaries}\label{subsec.prelim}
For more background on the theory of symmetric functions, we refer the reader to~\cite{Mac95}, and for more background on the representation theory of the symmetric group, we refer the reader to~\cite{fulton1997}.

A {\em partition} $\lambda$ of $n$, denoted by $\lambda \vdash n$, 
is a non-increasing sequence of positive integers $\lambda_1\geq\lambda_2\geq\cdots\geq \lambda_r>0$ whose sum is $n$.  Let $\calP$ denote the set of all partitions.
For $i\in \bbN$, let $m_i(\lambda)$ be the multiplicity of $i$ in $\lambda$, and let $z_\lambda = \prod_{i\geq1} i^{m_i(\lambda)} m_i(\lambda)!$.  

For $k\in \bbN$, the {\em $k$th power sum symmetric function} is
$$p_k=p_k(x_1,x_2,\ldots) = x_1^k + x_2^k + \cdots,$$
and the {\em power sum symmetric function} indexed by $\lambda$ is
$$p_\lambda= p_{\lambda_1}p_{\lambda_2}\cdots p_{\lambda_r}.$$
Let $\Lambda_\bbC$ denote the space of symmetric functions in the indeterminates $x_1, x_2, \ldots$.   Then $\{p_\lambda \mid \lambda \in \calP\}$ is a basis for $\Lambda_\bbC$.

Let $\fS_n$ denote the symmetric group on $n$ elements.  The irreducible representations of $\fS_n$ over $\bbC$ are indexed by the partitions of $n$, and are called {\em Specht modules}.  Let $\calS^\lambda$ denote the Specht module indexed by $\lambda$. 

The {\em Schur symmetric function} indexed by $\lambda$ satisfies
$$s_\lambda = \sum_\mu \chi^\lambda(\mu)z_\mu^{-1} p_\mu,$$
where $\chi^\lambda(\mu)$ is the character of the irreducible representation $\calS^\lambda$, evaluated on a permutation of cycle type $\mu$.  
The set $\{s_\lambda \mid \lambda \in \calP \}$ is another basis for $\Lambda_\bbC$.
By duality, 
one obtains the Murnaghan-Nakayama rule~\cite[I.7.8]{Mac95}, 
\begin{equation}\label{eqn.mnrule}
p_\mu = \sum_\lambda \chi^\lambda(\mu) s_\lambda,
\end{equation}
and combining this with the formula in equation~\eqref{thm.X_in_p} yields a formula for $X_G$ in terms of the Schur basis.  

\begin{remark}
If $\{u_\lambda\}$ is a basis for the space of symmetric functions, then $X_G = \sum_\lambda a_\lambda u_\lambda $ is said to be {\em $u$-positive} if $a_\lambda \geq 0$ for all partitions $\lambda$.  While it is clear that $X_G$ is not $p$-positive in general, it turns out that $X_G$ is $s$-positive for a certain class of graphs, although this is far from obvious following from the recipe for $X_G$ provided above.
~\cite{gasharov1996} obtained a remarkable combinatorial formula for $X_G$ in the Schur basis for the case when $G$ is an incomparability graph of a $(\mathbf{3}+\mathbf{1})$-free poset, and he used this to show that $X_G$ is $s$-positive in this case. 
\end{remark}

The {\em Grothendieck group} $R_n$ of representations of $\fS_n$ is the free abelian group on the isomorphism classes $[\calS^\lambda]$ of irreducible representations of $\fS_n$.  Let $R=\bigoplus_{n\geq0} R_n$.  If $[V]\in R_a$ and $[W]\in R_b$, then multiplication in $R$ is defined by
\begin{equation}
[V]\circ [W] = \left[\Ind_{\fS_a\times \fS_b}^{\fS_{a+b}} (V\otimes W)\right].
\end{equation}

Given irreducible modules $\calS^\lambda$ and $\calS^\mu$, the Littlewood-Richardson Rule~\cite[\S 5.2 Corollary 3]{fulton1997} is a combinatorial formula for computing their product
\begin{equation}\label{eqn.lrrule}
\left[\calS^\lambda\right]\circ \left[\calS^\mu\right] = \bigoplus_\nu c_{\lambda\mu}^\nu \left[\calS^\nu\right],
\end{equation}
where $c_{\lambda\mu}^\nu$ is the number of Littlewood-Richardson tableaux of shape $\nu/\lambda$ and weight $\mu$.  In other words, $c_{\lambda\mu}^\nu$ is the number of 
fillings of the skew shape $\nu/\lambda$ with $m_i(\mu)$ $i$'s, such that the filling is weakly increasing along rows, strictly increasing along columns, and the reading word obtained by reading the tableau
from top to bottom and from right to left is a lattice word (ie. in every initial subword, $i$ appears at least as often as $i+1$ for every $i$.) 

The morphism of graded rings given by sending the Specht modules to the Schur functions
\begin{equation}\label{eqn.frobchar}
\ch: R \rightarrow \Lambda_\bbC: [\calS^\lambda] \mapsto s_\lambda
\end{equation}
is an isomorphism~\cite[\S7.3 Theorem 1]{fulton1997}.

For $n\geq1$, it is a consequence of the Murnaghan-Nakayama rule~\eqref{eqn.mnrule} that
\begin{equation}\label{eqn.p_in_s}
\ch^{-1}(p_n) =  \sum_{i=0}^{n-1} (-1)^i\left[\calS^{(n-i, 1^i)}\right],
\end{equation}
is an alternating sum is over all hook-shaped partitions of $n$.

\section{The construction}\label{sec.construction}
\subsection{The state diagram}\label{subsec.statediagram}
Let $G$ be a graph with vertex set $V(G)= \{v_1,\ldots, v_n\}$ and edge set $E(G)= \{e_1,\ldots, e_m\}$.  

\begin{defn}
A {\em state} of $G$ is a spanning subgraph with a subset of edges $F\subseteq E$.  Let $|F|$ denote the number of edges in $F$.  
\end{defn}

A graph with $m$ edges therefore has $2^m$ states.  Let $Q(G)$ denote the set of all states of $G$.
We identify a state $F$ with the binary string $z_F=z_1\cdots z_n$, where 
$$z_i = \begin{cases} 0,&\hbox{if } e_i\in F,\\
1,&\hbox{if } e_i\notin F.
\end{cases}$$
With this identification, the states of $G$ form a Boolean lattice; that is, in the Hasse diagram of $Q(G)$, there is an edge from a state $F$ to a state $F'$ if and only if their associated binary strings differ in exactly one position.  
The Boolean lattice $Q(G)$ is graded by the number of edges in the states.

\begin{example} \label{eg.kthree_states}
Consider the $3$-cycle $G=K_3$ with vertex set $V=\{v_1,v_2,v_3\}$ and edge set $E=\{e_1,e_2,e_3\}$.
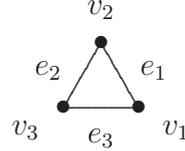
$$ 
\beginpicture
	\setcoordinatesystem units <1cm,1cm>         
	\setplotarea x from 0 to 0, y from 0 to 0    
	\multiput{$\bullet$} at 0 0 1 0 0.5 0.867 /
	\put{$v_1$} at 1.5 -.3 \put{$v_2$} at .5 1.3 \put{$v_3$} at -.5 -.3 
	\put{$e_1$}[c] at 1.2 .5 \put{$e_2$}[c] at -.2 .5 
	\put{$e_3$}[c] at .5 -.4
	\plot 0 0 1 0 0.5 .867 0 0 /
\endpicture$$

The eight states of $K_3$ and their associated binary strings are shown below, together with the signed per-edge maps (to be defined in Section~\ref{subsec.peredge}):
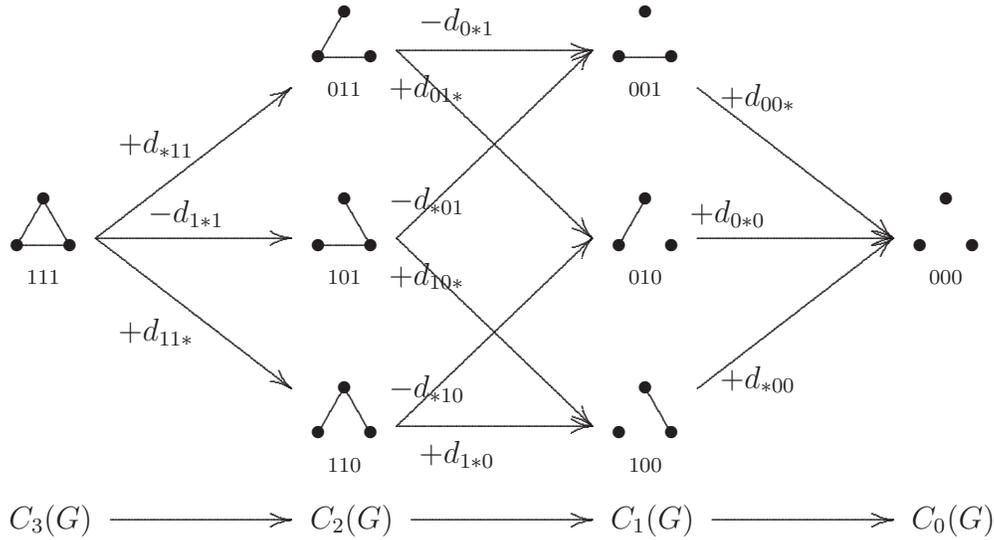
$$ 
\beginpicture
	\setcoordinatesystem units <4cm,2.5cm>         
	\setplotarea x from 0 to 0, y from -1.75 to 1.5    

\arrow <8pt> [.2,.67] from 0.1 0 to .75 .8 	\put{$+d_{*11}$} at .3 .5 
\arrow <8pt> [.2,.67] from 0.1 0 to .75 0    \put{$-d_{1*1}$} at .4 .13
\arrow <8pt> [.2,.67] from 0.1 0 to .75 -.8  \put{$+d_{11*}$} at .3 -.5

\arrow <8pt> [.2,.67] from 1.1 1 to 1.75 1   \put{$-d_{0*1}$}[c] at 1.3 1.15
\arrow <8pt> [.2,.67] from 1.1 1 to 1.75 0   \put{$+d_{01*}$}[c] at 1.2 .8
\arrow <8pt> [.2,.67] from 1.1 0 to 1.75 1   \put{$-d_{*01}$}[c] at 1.2 .2
\arrow <8pt> [.2,.67] from 1.1 0 to 1.75 -1  \put{$+d_{10*}$}[c] at 1.2 -.2
\arrow <8pt> [.2,.67] from 1.1 -1 to 1.75 0  \put{$-d_{*10}$}[c] at 1.2 -.8
\arrow <8pt> [.2,.67] from 1.1 -1 to 1.75 -1 \put{$+d_{1*0}$}[c] at 1.3 -1.15

\arrow <8pt> [.2,.67] from 2.1 .8 to 2.75 0  \put{$+d_{00*}$}[c] at 2.3 .75
\arrow <8pt> [.2,.67] from 2.1 0 to 2.75 0   \put{$+d_{0*0}$}[c] at 2.2 .13
\arrow <8pt> [.2,.67] from 2.1 -.8 to 2.75 0 \put{$+d_{*00}$}[c] at 2.3 -.75 

\put{$\beginpicture
	\setcoordinatesystem units <.7cm,.7cm>         
	\setplotarea x from 0 to 2, y from 0 to 0    
	\multiput{$\bullet$} at 0 0 1 0 0.5 0.867 /
	\plot 0 0 1 0 0.5 .867 0 0 /
	\put{\tiny111} at .5 -.6
	\endpicture$} at 0 0

\put{$\beginpicture
	\setcoordinatesystem units <.7cm,.7cm>         
	\setplotarea x from 0 to 2, y from 0 to 0    
	\multiput{$\bullet$} at 0 0 1 0 0.5 0.867 /
	\plot 1 0 0 0 0.5 .867 /
	\put{\tiny011} at .5 -.6
	\endpicture$} at 1 1
\put{$\beginpicture
	\setcoordinatesystem units <.7cm,.7cm>         
	\setplotarea x from 0 to 2, y from 0 to 0    
	\multiput{$\bullet$} at 0 0 1 0 0.5 0.867 /
	\plot 0.5 .867 1 0 0 0 /
	\put{\tiny101} at .5 -.6
	\endpicture$} at 1 0	
\put{$\beginpicture
	\setcoordinatesystem units <.7cm,.7cm>         
	\setplotarea x from 0 to 2, y from 0 to 0    
	\multiput{$\bullet$} at 0 0 1 0 0.5 0.867 /
	\plot 0 0 0.5 .867 1 0 /
	\put{\tiny110} at .5 -.6
	\endpicture$} at 1 -1

\put{$\beginpicture
	\setcoordinatesystem units <.7cm,.7cm>         
	\setplotarea x from 0 to 2, y from 0 to 0    
	\multiput{$\bullet$} at 0 0 1 0 0.5 0.867 /
	\plot 0 0 1 0 /
	\put{\tiny001} at .5 -.6
	\endpicture$} at 2 1	
\put{$\beginpicture
	\setcoordinatesystem units <.7cm,.7cm>         
	\setplotarea x from 0 to 2, y from 0 to 0    
	\multiput{$\bullet$} at 0 0 1 0 0.5 0.867 /
	\plot 0 0 0.5 .867 /
	\put{\tiny010} at .5 -.6
	\endpicture$} at 2 0
\put{$\beginpicture
	\setcoordinatesystem units <.7cm,.7cm>         
	\setplotarea x from 0 to 2, y from 0 to 0    
	\multiput{$\bullet$} at 0 0 1 0 0.5 0.867 /
	\plot 1 0 0.5 .867 /
	\put{\tiny100} at .5 -.6
	\endpicture$} at 2 -1

\put{$\beginpicture
	\setcoordinatesystem units <.7cm,.7cm>         
	\setplotarea x from 0 to 2, y from 0 to 0    
	\multiput{$\bullet$} at 0 0 1 0 0.5 0.867 /
	\put{\tiny000} at .5 -.6
	\endpicture$} at 3 0
\put{$C_3(G)$}[c] at -.05 -1.5
\put{$C_2(G)$}[c] at .95 -1.5
\put{$C_1(G)$}[c] at 1.95 -1.5
\put{$C_0(G)$}[c] at 2.95 -1.5
\arrow <8pt> [.2,.67] from 0.15 -1.5 to 0.75 -1.5 
\arrow <8pt> [.2,.67] from 1.15 -1.5 to 1.75 -1.5 	
\arrow <8pt> [.2,.67] from 2.15 -1.5 to 2.75 -1.5 	
\endpicture
$$
\end{example}

\subsection{The graded modules}
We shall assign a graded $\fS_n$-module to each state.  Let $\calS^\lambda$ denote the irreducible $\fS_n$-module indexed by the partition $\lambda\vdash n$.  For $a=1$, let $\calL_a = \calS^{(a)}$, and for $a\in \bbZ_{\geq2}$, let $\calL_a$ denote the graded $\fS_n$-module
\begin{equation}
\calL_a = \bigoplus_{i=0}^{a-1} \calS^{(a-i,1^i)}.
\end{equation} 
Since $\largewedge^i \calS^{(a-1,1)} \cong \calS^{(a-i, 1^i)}$, we may also think of $\calL_a$ as the exterior algebra $\largewedge^* \calS^{(a-1,1)}$.

Suppose $F\subseteq E(G)$ is a state with $r$ connected components $B_1,\ldots, B_r$ of sizes $b_1, \ldots, b_r$ respectively.   To $F$, we assign the graded $\fS_n$-module 
\begin{equation}
\calM_F = \Ind_{\fS_{B_1}\times \cdots \times \fS_{B_r}}^{\fS_{V(G)}} 
\left(\calL_{b_1}\otimes \cdots \otimes \calL_{b_r}\right),
\end{equation}
where $\fS_{B_1}\times \cdots \times \fS_{B_r}$ is a Young subgroup of $\fS_{V(G)}\cong \fS_n$.  

For the sake of simpler notation, we will sometimes use the shorthand
\begin{equation}
\Ind_{B_1|\cdots |B_r} = \Ind_{\fS_{B_1}\times \cdots \times \fS_{B_r}}^{\fS_{V(G)}}.
\end{equation}

\begin{defn}
For $i\geq0$, the $i$th {\em chain module} for $G$ is
$$C_i(G) = \bigoplus_{|F|=i} \calM_F. $$
More precisely, since $\calM_F = \bigoplus_{j\geq0} \left(\calM_F\right)_j$ is graded, then for $i,j\geq0$, we can define
$$C_{i,j}(G) = \bigoplus_{|F|=i} \left(\calM_F\right)_j.$$
Observe that $C_{i,j}(G)=0$ for all $j>i$.
\end{defn}

\begin{example}\label{eg.kthree_objects}
Continuing the $K_3$ example, the graded $\fS_3$-modules corresponding to the states are:
$$ 
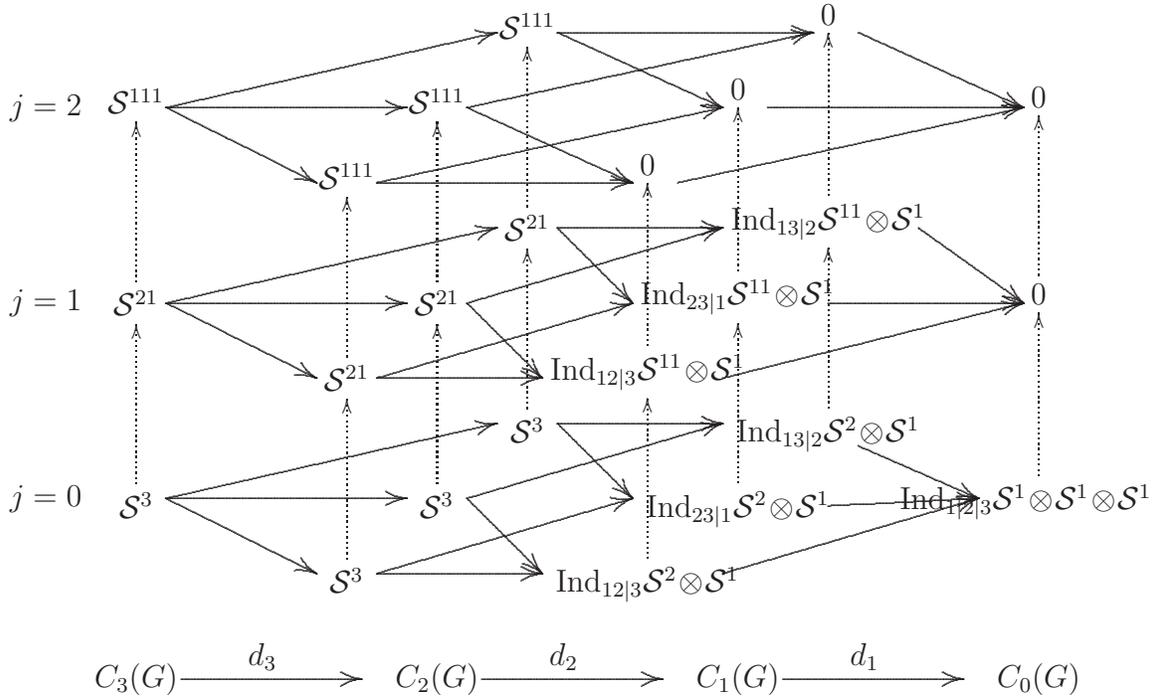
\beginpicture
	\setcoordinatesystem units <4cm,2cm>         
	\setplotarea x from -.45 to 2.4, y from -3 to 2    

\arrow <8pt> [.2,.67] from .1 1.3 to 1.2 1.8 
\arrow <8pt> [.2,.67] from .1 1.3 to .9 1.3 
\arrow <8pt> [.2,.67] from .1 1.3 to .6 .8 

\arrow <8pt> [.2,.67] from 1.4 1.8 to 2.25 1.8 
\arrow <8pt> [.2,.67] from 1.4 1.8 to 1.95 1.3 
\arrow <8pt> [.2,.67] from 1.1 1.3 to 2.25 1.8 
\arrow <8pt> [.2,.67] from 1.1 1.3 to 1.65 .8 
\arrow <8pt> [.2,.67] from 0.8 .8 to 1.95 1.3 
\arrow <8pt> [.2,.67] from 0.8 .8 to 1.65 .8 

\arrow <8pt> [.2,.67] from 2.4 1.8 to 2.95 1.3 
\arrow <8pt> [.2,.67] from 2.1 1.3 to 2.95 1.3 
\arrow <8pt> [.2,.67] from 1.8 .8 to 2.95 1.3 

\arrow <8pt> [.2,.67] from .1 0 to 1.2 .5 
\arrow <8pt> [.2,.67] from .1 0 to .9 0 
\arrow <8pt> [.2,.67] from .1 0 to .6 -.5 

\arrow <8pt> [.2,.67] from 1.4 .5 to 1.95 .5 
\arrow <8pt> [.2,.67] from 1.4 .5 to 1.65 0 
\arrow <8pt> [.2,.67] from 1.1 0 to 1.95 .5 
\arrow <8pt> [.2,.67] from 1.1 0 to 1.35 -.5 
\arrow <8pt> [.2,.67] from 0.8 -.5 to 1.65 0 
\arrow <8pt> [.2,.67] from 0.8 -.5 to 1.35 -.5 

\arrow <8pt> [.2,.67] from 2.6 .5 to 2.95 0 
\arrow <8pt> [.2,.67] from 2.3 0 to 2.95 0 
\arrow <8pt> [.2,.67] from 1.95 -.5 to 2.95 0 

\arrow <8pt> [.2,.67] from .1 -1.3 to 1.2 -.8 
\arrow <8pt> [.2,.67] from .1 -1.3 to .9 -1.3 
\arrow <8pt> [.2,.67] from .1 -1.3 to .6 -1.8 

\arrow <8pt> [.2,.67] from 1.4 -.8 to 1.95 -.8 
\arrow <8pt> [.2,.67] from 1.4 -.8 to 1.65 -1.3 
\arrow <8pt> [.2,.67] from 1.1 -1.3 to 1.95 -.8 
\arrow <8pt> [.2,.67] from 1.1 -1.3 to 1.35 -1.8 
\arrow <8pt> [.2,.67] from 0.8 -1.8 to 1.65 -1.3 
\arrow <8pt> [.2,.67] from 0.8 -1.8 to 1.35 -1.8 

\arrow <8pt> [.2,.67] from 2.4 -.95 to 2.8 -1.3 
\arrow <8pt> [.2,.67] from 2.3 -1.35 to 2.8 -1.3 
\arrow <8pt> [.2,.67] from 1.95 -1.8 to 2.8 -1.3 

\put{\xymatrixrowsep{60pt}\xymatrix{
\calS^{111}\\
\calS^{21} \ar@{.>}[u]\\
\calS^{3} \ar@{.>}[u]
}} at 0 0

\put{\xymatrixrowsep{60pt}\xymatrix{
\calS^{111}\\
\calS^{21} \ar@{.>}[u]\\
\calS^{3} \ar@{.>}[u]
}} at 1.3 .5
\put{\xymatrixrowsep{60pt}\xymatrix{
\calS^{111}\\
\calS^{21} \ar@{.>}[u]\\
\calS^{3} \ar@{.>}[u]
}} at 1 0
\put{\xymatrixrowsep{60pt}\xymatrix{
\calS^{111}\\
\calS^{21} \ar@{.>}[u]\\
\calS^{3} \ar@{.>}[u]
}} at .7 -.5

\put{\xymatrixrowsep{60pt}\xymatrix{
0\\
\Ind_{13|2} \calS^{11}\!\otimes\! \calS^{1} \ar@{.>}[u]\\
\Ind_{13|2} \calS^{2}\!\otimes\! \calS^{1} \ar@{.>}[u]
}} at 2.3 .5
\put{\xymatrixrowsep{60pt}\xymatrix{
0\\
\Ind_{23|1} \calS^{11}\!\otimes\!\calS^{1} \ar@{.>}[u]\\
\Ind_{23|1} \calS^{2}\!\otimes\! \calS^{1} \ar@{.>}[u]
}} at 2 0
\put{\xymatrixrowsep{60pt}\xymatrix{
0\\
\Ind_{12|3} \calS^{11}\!\otimes\! \calS^{1} \ar@{.>}[u]\\
\Ind_{12|3} \calS^{2}\!\otimes\! \calS^{1} \ar@{.>}[u]
}} at 1.7 -.5

\put{\xymatrixrowsep{60pt}\xymatrix{
0\\
0 \ar@{.>}[u]\\
\!\!\!\!\!
\Ind_{1|2|3} \calS^{1}\!\otimes\! \calS^{1}\!\otimes\!\calS^{1} \ar@{.>}[u]
}} at 3 0
	
\put{$j=2$}[c] at -.3 1.3
\put{$j=1$}[c] at -.3 0
\put{$j=0$}[c] at -.3 -1.3

\put{$C_3(G)$}[c] at 0 -2.5
\put{$C_2(G)$}[c] at 1 -2.5
\put{$C_1(G)$}[c] at 2 -2.5
\put{$C_0(G)$}[c] at 3 -2.5
\arrow <8pt> [.2,.67] from 0.15 -2.5 to 0.75 -2.5 
	\put{$d_{3}$} at .42 -2.35
\arrow <8pt> [.2,.67] from 1.15 -2.5 to 1.75 -2.5 	
	\put{$d_{2}$} at 1.42 -2.35
\arrow <8pt> [.2,.67] from 2.15 -2.5 to 2.75 -2.5 	
	\put{$d_{1}$} at 2.42 -2.35
\endpicture
$$

\end{example}

\subsection{The per-edge maps}\label{subsec.peredge}
So far, the modules $\calM_F$ are defined up to isomorphism.  In the following, we show that for each edge $\varepsilon$ in the Hasse diagram of $Q(G)$ from a state $F$ to a state $F'$, there is an obvious choice of a graded $\fS_n$-module morphism $d_\varepsilon: \calM_F \rightarrow \calM_{F'}$.  These are the {\em per-edge maps}.  

There are two cases to consider.
Suppose $F' = F-e$ where $e\in E(G)$.

\subsubsection*{\bf Case 1}
The edge $e$ is incident to vertices in the same connected component of $F'$.

Since $\calM_F$ and $\calM_{F'}$ are canonically isomorphic,
we define $d_\varepsilon:\calM_F\rightarrow \calM_{F'} $ to be the identity map.

\subsubsection*{\bf Case 2}
The edge $e$ is incident to vertices in different connected components of $F'$.

First, consider the simplest case where $F$ consists of one connected component, and $F'$ consists of two connected components $A$ and $B$.
Suppose $|A|=a$ and $|B|=b$, so that $a+b=n$. 
Given that $\calM_F \cong \calL_n$ and $\calM_{F'} \cong \Ind_{\fS_a\times\fS_b} (\calL_a \otimes \calL_b)$,  
we shall show that there is a canonical choice of an element in
$\Hom_{\fS_n}\!\left(\calM_F, \calM_{F'}\right)$.

\begin{lemma} \label{lem.res}
Let $\mathbf{1}_W$ denote the trivial representation of the group $W$.  Then
$$\Res_{\fS_a\times \fS_b}^{\fS_n} \calS^{(n-1,1)}
\cong 
	\left(\left(\calS^{(a-1,1)} \otimes \mathbf{1}_{\fS_b}\right) 
	\oplus \left(\mathbf{1}_{\fS_a} \otimes \calS^{(b-1,1)}\right)\right) 
	\oplus \mathbf{1}_{\fS_a \times \fS_b}. 
$$
\end{lemma}

\proof
The permutation representation of $\fS_n$ is the module
$$\calM^{(n-1,1)} \cong \Ind_{\fS_{n-1}}^{\fS_n} \mathbf{1}_{\fS_{n-1}}.$$
Since $\calM^{(n-1,1)} \cong  \calS^{(n-1,1)}\oplus \calS^{(n)}$ and
$\calS^{(n)}\cong \mathbf{1}_{\fS_n}$,
then
$$\Res_{\fS_a\times \fS_b}^{\fS_n}\calM^{(n-1,1)}
\cong \left(\Res_{\fS_a\times \fS_b}^{\fS_n} \calS^{(n-1,1)}\right)
	\oplus \mathbf{1}_{\fS_a\times \fS_b}.
$$
On the other hand,
\begin{align*}
\Res_{\fS_a\times \fS_b}^{\fS_n}\calM^{(n-1,1)}
&\cong 
\left(\calM^{(a-1,1)} \otimes \mathbf{1}_{\fS_b}\right) 
	\oplus \left(\mathbf{1}_{\fS_a} \otimes \calM^{(b-1,1)}\right) \\
&\cong 
\left(\calS^{(a-1,1)} \otimes \mathbf{1}_{\fS_b}\right) 
	\oplus \left(\mathbf{1}_{\fS_a} \otimes \calS^{(b-1,1)}\right)
	\oplus \left(\mathbf{1}_{\fS_a \times \fS_b}\right)^{\oplus2}.
\end{align*}
Comparing the two expressions for $\Res_{\fS_a\times \fS_b}^{\fS_n}\calM^{(n-1,1)}$ yields the result.
\qed

\begin{lemma}  \label{lem.frobrec}
Let 
$\calT = \left(\calS^{(a-1,1)} \otimes \mathbf{1}_{\fS_b}\right) 
	\oplus \left(\mathbf{1}_{\fS_a} \otimes \calS^{(b-1,1)}\right).$
Then
$$\Hom_{\fS_n}\!\left( \calM_F, \calM_{F'}\right)	
\cong \Hom_{\fS_a\times \fS_b}\!
	\left(\largewedge^*\calT \oplus (\largewedge^*\calT)[1], \largewedge^*\calT
	\right),
$$
where $(\largewedge^* \calT)[1]$ denotes the first translate of $\largewedge^* \calT$.
\end{lemma}
\proof
By Frobenius Reciprocity, 
\begin{align*}
\Hom_{\fS_n}\!\left( \calM_F, \calM_{F'}\right)	
&= \Hom_{\fS_n}\!
	\left( \calL_n,
	\Ind_{\fS_a\times \fS_b}^{\fS_n} (\calL_a\otimes \calL_b)
	\right)\\
&= \Hom_{\fS_a\times \fS_b}\!\left( 
	\Res_{\fS_a\times \fS_b}^{\fS_n}(\calL_n),
	\calL_a\otimes \calL_b
	\right)
\end{align*}
is a natural isomorphism of vector spaces.  Note that
$$\calL_a\otimes \calL_b
=\largewedge^*\calS^{(a-1,1)} \otimes \largewedge^*\calS^{(b-1,1)}
=\largewedge^*(\calS^{(a-1,1)}\oplus \calS^{(b-1,1)})
\cong 
\largewedge^*\calT.  $$
Also, by Lemma~\ref{lem.res},
\begin{align*}
\Res_{\fS_a\times \fS_b}^{\fS_n} \calL_n
	&= \Res_{\fS_a\times \fS_b}^{\fS_n}
	\left(\largewedge^* \calS^{(n-1,1)}\right)
	= \largewedge^*\!
	\left(\Res_{\fS_a\times \fS_b}^{\fS_n} \calS^{(n-1,1)}\right)\\
	&= \largewedge^*\!
	\left( \calT \oplus \mathbf{1}_{\fS_a\times \fS_b}\right)
	=\largewedge^* \calT  
	\otimes \largewedge^* \mathbf{1}_{\fS_a\times \fS_b}\\
	&=\largewedge^* \calT \otimes \left(\mathbf{1}_{\fS_a\times \fS_b} 
		\oplus \mathbf{1}_{\fS_a\times \fS_b}[1] \right)\\
	&=\largewedge^* \calT \oplus (\largewedge^* \calT)[1].
\end{align*}
From this, the result follows.
\qed

Since 
$\Hom_{\fS_n}\!\left( \calM_F, \calM_{F'}\right)	
=\Hom_{\fS_a\times \fS_b}\!
	\left(\largewedge^*\calT \oplus (\largewedge^*\calT)[1],
	\largewedge^*\calT
	\right),
$
we choose the element $d_\varepsilon\in \Hom_{\fS_n}\!\left( \calM_F,\calM_{F'}\right)$ to be the map that corresponds to the $(\fS_a\times\fS_b)$-module map that is identity on $\largewedge^*\calT$, and zero on $(\largewedge^*\calT)[1]$.

It remains to consider the general case when $F$ has more than one connected component.  The definition of the per-edge map in this case is achieved through induction on the two-component case.

Suppose $F$ is a state with $r$ connected components $B_1, \ldots, B_r$ of sizes $b_1,\ldots, b_r$.  Further suppose that the removal of the edge $e\in E(G)$ decomposes $B_r$ into two components $A$ and $B$ of sizes $a$ and $b$ respectively.

Let $d_\zeta: \calL_{a+b} \rightarrow \Ind_{\fS_A\times\fS_B}^{\fS_{B_r}} \left(\calL_{a}\otimes \calL_b\right)$ be the per-edge map defined previously, and
let $\calN = \calL_{b_1}\otimes \cdots \otimes \calL_{b_{r-1}}$.
The map $d_\varepsilon:\calM_F\rightarrow \calM_{F'}$ is chosen to be
\begin{equation}\label{eqn.indd}
d_\varepsilon 
=\Ind_{\fS_{B_1}\times\cdots\times\fS_{B_{r-1}}\times\fS_{B_r}}^{\fS_{V(G)}}
\left(\mathrm{id}_\calN \otimes d_\zeta\right).
\end{equation}

\subsubsection{Signed per-edge maps}
Lastly, we define signed per-edge maps. 
Suppose the associated binary strings $z_F$ and $z_{F'}$ differ in the $j$th position.  Identify $\varepsilon$ with the sequence $\varepsilon_1\cdots \varepsilon_m$ defined by $\varepsilon_j = *$, and
$\varepsilon_i = (z_F)_i = (z_{F'})_i$ for $i\neq j$.  See Example~\ref{eg.kthree_states}.

\begin{defn}
If an edge $\varepsilon=\varepsilon_1\cdots \varepsilon_m$ in the Hasse diagram of $Q(G)$ has $k$ ones after the $*$ entry, then the {\em signed per-edge map} that corresponds to $\varepsilon$ is
$$\mathrm{sgn}(\varepsilon)d_\varepsilon = (-1)^k d_\varepsilon.
$$
\end{defn}

\begin{remark}
This choice of a sign convention on the edges of $Q(G)$ guarantees that any induced sublattice isomorphic to a square will have three edges of the same sign and one edge of the opposite sign, which is crucial in the proof that $d$ is indeed a differential map, as will be seen in the next section.
\end{remark}

\subsection{The chain complex of graded modules}\label{subsec.chaincx}
Now that we have assigned a graded $\fS_n$-module $\calM_F$ to each state $F$ of $G$, and defined signed per-edge maps, we can make the following definition.

\begin{defn}
For $i\geq0$, define $d_i:C_i(G) \rightarrow C_{i-1}(G)$ by letting
$$d_i = \sum_{\varepsilon} \mathrm{sgn}(\varepsilon)d_\varepsilon,$$
where the sum is over all edges $\varepsilon$ in the Hasse diagram of $Q(G)$ joining a state with $i$ edges to a state with $i-1$ edges.   We also define $d_{i,j}:C_{i,j}(G) \rightarrow C_{i-1,j}(G)$ to be the map $d_i$ in the $j$th grading.
\end{defn}

\begin{prop}\label{prop.differential}
This defines a differential; that is, $d^2=0$.
\end{prop}

\proof Since $d$ is defined via per-edge maps, it suffices to show that $d^2=0$ on the per-edge maps.  This is achieved through case checking; up to signs, the three cases that must be considered depend on how the removal of two edges from a state of $G$ will disconnect graph components.

Let $F$ be a state of $G$, and suppose $e_1$ and $e_2$ are the two edges to be removed from $F$.  We use the notation $F_i = F-e_i$, and $F_{12} = F-e_1-e_2$.  Label the four edges in the Hasse diagram of $Q(G)$ as in the following diagram.  It is straightforward to verify that the sign convention on the per-edge maps forces each `square' in the Hasse diagram $Q(G)$ to have three of one sign, and one of the opposite sign.

We focus on the case with the signs indicated in the following diagram; the other three cases are completely analogous.

$$
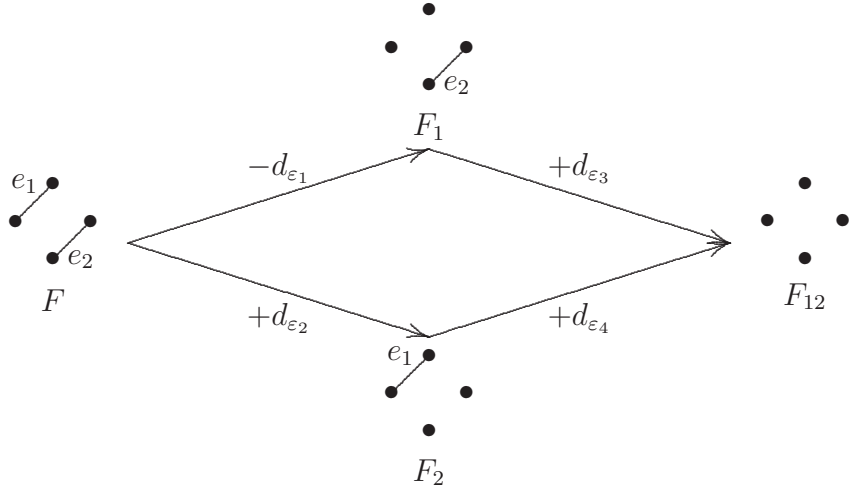
\beginpicture
	\setcoordinatesystem units <4cm,1.25cm>         
	\setplotarea x from 0 to 2, y from -.5 to 2    
\put{$
\beginpicture
	\setcoordinatesystem units <.5cm,.5cm>
	\setplotarea x from -1 to 3, y from -1 to 2  
	\put{$\bullet$} at 0 1 \put{$\bullet$} at 1 0 \put{$\bullet$} at 1 2
	\put{$\bullet$} at 2 1
	\plot 0 1 1 2 / \plot 1 0 2 1 /
	\put{$e_1$} at .25 2 \put{$e_2$} at 1.75 0
	\put{$F$}[c] at 1 -1.1
\endpicture
$}[r] at 0 1

\put{$
\beginpicture
	\setcoordinatesystem units <.5cm,.5cm>
	\setplotarea x from -1 to 3, y from -1 to 2  
	\put{$\bullet$} at 0 1 \put{$\bullet$} at 1 0 \put{$\bullet$} at 1 2
	\put{$\bullet$} at 2 1
	\plot 1 0 2 1 / \put{$e_2$} at 1.75 0
	\put{$F_1$}[c] at 1 -1.1
\endpicture
$}[b] at 1 2.1

\put{$
\beginpicture
	\setcoordinatesystem units <.5cm,.5cm>
	\setplotarea x from -1 to 3, y from -1 to 2 
	\put{$\bullet$} at 0 1 \put{$\bullet$} at 1 0 \put{$\bullet$} at 1 2
	\put{$\bullet$} at 2 1
	\plot 0 1 1 2 / \put{$e_1$} at .25 2
	\put{$F_2$}[c] at 1 -1.1
\endpicture
$}[ct] at 1 -0.1

\put{$
\beginpicture
	\setcoordinatesystem units <.5cm,.5cm>
	\setplotarea x from -1 to 3, y from -1 to 2 
	\put{$\bullet$} at 0 1 \put{$\bullet$} at 1 0 \put{$\bullet$} at 1 2
	\put{$\bullet$} at 2 1
	\put{$F_{12}$}[c] at 1 -1
\endpicture
$}[l] at 2 1

\arrow <8pt> [.2,.67] from 0 1  to 1 2 \put{$-d_{\varepsilon_1}$} at .5 1.8
\arrow <8pt> [.2,.67] from 0 1  to 1 0 \put{$+d_{\varepsilon_2}$} at .5 .2
\arrow <8pt> [.2,.67] from 1 2  to 2 1 \put{$+d_{\varepsilon_3}$} at 1.5 1.8
\arrow <8pt> [.2,.67] from 1 0  to 2 1 \put{$+d_{\varepsilon_4}$} at 1.5 .2
\endpicture
$$

Now we consider the three cases which depend on how the removal of $e_1$ and $e_2$ disconnects graph components.
 Let $\kappa(F)$ denote the number of connected components of $F$.  Observe that $\kappa(F_{12}) \leq \kappa(F) +2$.  In each case, we verify that $d_{\varepsilon_4}\circ d_{\varepsilon_2} - d_{\varepsilon_3} \circ d_{\varepsilon_1} =0$.

\noindent {\bf Case 1.} 
$\kappa(F_{12}) = \kappa(F)$.  
In this case, the modules corresponding to all four states are equal, and each per-edge map is the identity map (up to signs), so for any $x\in \calM_F$, 
$$x\mapsto (-x,x) \mapsto x-x=0. $$

\noindent {\bf Case 2.} 
$\kappa(F_{12}) = \kappa(F)+1$.  There are two subcases to consider.
First, suppose removing $e_1$ disconnects a graph component, but removing $e_2$ does not.
Since the removal of $e_2$ does not disconnect components, then $\calM_F = \calM_{F_2}$, $\calM_{F_1} = \calM_{F_{12}}$, and 
the two per-edge maps $d_{\varepsilon_2}$ and $d_{\varepsilon_3}$ are the identity maps.
If the removal of $e_1$ disconnects a component of size $n$ into components of sizes $a$ and $b$, then the maps $d_{\varepsilon_1}$ and $d_{\varepsilon_4}$ are each induced from the same map $ \calL_n \rightarrow \Ind_{\fS_a\times\fS_b}^{\fS_n} \left(\calL_a\otimes \calL_b\right)$.  Therefore, $d_{\varepsilon_1}=d_{\varepsilon_4}$, 
and 
$$d_{\varepsilon_4} \circ \mathrm{id}_{\calM_F}
- \mathrm{id}_{\calM_{F_1}}\circ d_{\varepsilon_1} 
=0.$$

Second, suppose that removing each of $e_1$ or $e_2$ does not disconnect a graph component, but that removing both edges simultaneously disconnects a component.  Then  $\calM_F = \calM_{F_1}=\calM_{F_2}$, and the two per-edge maps $d_{\varepsilon_1}$ and $d_{\varepsilon_2}$ are the identity maps up to sign.
If the removal of $e_1$ and $e_2$ disconnects a component of size $n$ into components of sizes $a$ and $b$, then the maps $d_{\varepsilon_3}$ and $d_{\varepsilon_4}$ are each induced from the same map $ \calL_n \rightarrow \Ind_{\fS_a\times\fS_b}^{\fS_n} \left(\calL_a\otimes \calL_b\right)$.
Therefore, $d_{\varepsilon_3}=d_{\varepsilon_4}$, 
and 
$$d_{\varepsilon_4} \circ id_{\calM_F}  
- d_{\varepsilon_3}\circ id_{\calM_F}
= 0.$$

\noindent {\bf Case 3.} $\kappa(F_{12}) = \kappa(F)+2$.  There are two subcases to consider.
First, suppose $e_1$ and $e_2$ are in different components $A$ and $B$ of $F$, respectively.  Let $|A|=a$ and $|B|=b$.
For simplicity, we assume that $F$ only has two components; the general case is handled by induction~\eqref{eqn.indd}.  
So $\calM_A = \calL_a$, $\calM_B = \calL_b$, and $\calM_F = \Ind_{\fS_a \times \fS_b}^{\fS_{a+b}} (\calL_a \otimes \calL_b)$.  Let $f:\calM_A \rightarrow \calM_{A-e_1}$ and $g:\calM_{B}\rightarrow \calM_{B-e_2}$.
We then have the following maps of $(\fS_a\times\fS_b)$-modules: 
$$\xymatrix{
&\calM_{A-e_1}\otimes\calM_B
	\ar[dr]^-{+\mathrm{id}\otimes g}&\\
\calM_A\otimes \calM_B
	\ar[ur]^-{-f\otimes \mathrm{id}}
	\ar[dr]_-{+\mathrm{id}\otimes g}
&&\calM_{A-e_1}\otimes\calM_{B-e_2}\\
&\calM_{A}\otimes\calM_{B-e_2}
	\ar[ur]_-{+f\otimes \mathrm{id}}&
}$$
so that 
$$(f\circ \mathrm{id}_{\calM_A})\otimes(\mathrm{id}_{\calM_{B-e_2}}\circ g) 
- (\mathrm{id}_{\calM_{A-e_1}}\circ f)\otimes(g\circ \mathrm{id}_{\calM_{B}}) 
= 0.$$
Induce this to $\fS_{a+b}$-modules to conclude that $d_{\varepsilon_4}\circ d_{\varepsilon_2} - d_{\varepsilon_3}\circ d_{\varepsilon_1}=0$.

Lastly, suppose $e_1$ and $e_2$ are in the same component of $F$.  Again for simplicity, we assume that $F$ consists of only one component, and the general case is handled by induction.  Suppose $F_{12}$ consists of the three components $A_1, A_2,$ and $A_3$ of sizes $a,b$ and $c$ respectively, while $F_1$ consists of the components $A_1$ and $A_2\cup A_3$ and $F_2$ consists of the components $A_1\cup A_2$ and $A_3$.
Let $n=a+b+c$.  We have the following maps of $\fS_n$-modules:
$$\xymatrix{
&\calM_{F_1}=\Ind \left(\calL_{a}\otimes \calL_{b+c}\right)
	 \ar[dr]^-{+d_{\varepsilon_3}}&\\
\calM_F=\calL_n
	\ar[ur]^-{-d_{\varepsilon_1}}
	\ar[dr]_-{+d_{\varepsilon_2}}
&&
\calM_{F_{12}}
	=\Ind \left(\calL_{a}\otimes \calL_{b}\otimes \calL_{c} \right)\\
&\calM_{F_2}=\Ind \left( \calL_{a+b}\otimes \calL_{c}\right) 
	\ar[ur]_-{+d_{\varepsilon_4}}&
}$$
Recall from equation~\ref{eqn.indd} that the map $d_{\varepsilon_3}$ is induced from $\mathrm{id}_{\calL_{a}} \otimes (\calL_{b+c}\rightarrow \Ind (\calL_{b}\otimes \calL_{c}))$,
so by Frobenius reciprocity, the $\fS_n$-module map $d_{\varepsilon_3}\circ d_{\varepsilon_1}$ corresponds to the $(\fS_{a}\times \fS_{b+c})$-module map
$$\xymatrix{
\Res^{\fS_n}_{\fS_{a}\times \fS_{b+c}} (\calL_n) \ar[r] & 
\calL_{a}\otimes \calL_{b+c} \ar[r] &
\calL_{a}\otimes \Ind_{\fS_{b}\times \fS_{c}}^{\fS_{b+c}} (\calL_{b}\otimes \calL_{c}).
}$$
Applying Frobenius reciprocity one more time, this map then in turn corresponds to the 
$(\calL_{a}\times \calL_{b}\times \calL_{c})$-module map
$$\xymatrix{
\psi:\Res_{\fS_{a}\times \fS_{b}\times \fS_{c}}^{\fS_n} (\calL_n) \ar[r] & 
\calL_{a}\otimes \Res^{\fS_{b+c}}_{\fS_{b}\times\fS_{c}}(\calL_{b+c}) \ar[r] &
\calL_{a}\otimes \calL_{b}\otimes \calL_{c}.
}$$
Letting
$\calL_{b}\otimes \calL_{c} = \largewedge^* \calU$, where
$\calU =
\left( \calS^{(b-1,1)}\otimes \mathbf{1}_{\fS_{c}}\right)
\oplus 
\left( \mathbf{1}_{\fS_{b}}\otimes \calS^{(c-1,1)}\right)$, 
then
$$\calL_{a}\otimes \Res^{\fS_{b+c}}_{\fS_{b}\times\fS_{c}} (\calL_{b+c})
=\calL_{a}\otimes(\largewedge^* \calU \oplus \largewedge^* \calU[1]).
$$
Also, letting $\calL_{a}\otimes \calL_{b+c} = \largewedge^* \calT$, where
$\calT  = \left( \calS^{(a-1,1)}\otimes \mathbf{1}_{\fS_{b+c}}\right)
\oplus 
\left( \mathbf{1}_{\fS_{a}}\otimes \calS^{(b+c-1,1)}\right)$, 
then
$$\Res_{\fS_{a}\times \fS_{b}\times \fS_{c}}^{\fS_n} (\calL_n)
= (\calL_{a}\otimes\largewedge^* \calU) 
\oplus \left(\calL_{a}\otimes\largewedge^* \calU\right)^{\oplus2} [1]
\oplus (\calL_{a}\otimes\largewedge^* \calU)[2],$$
therefore, the composite map $\psi:\Res(\calL_n)\rightarrow \calL_{a}\otimes \calL_{b}\otimes\calL_{c}$ is the identity map on $\calL_{a}\otimes\largewedge^*\calU$, and is the zero map on the remaining factors. 

Similarly, the $\fS_n$-module map $d_{\varepsilon_4}\circ d_{\varepsilon_2}$ corresponds to the 
$(\calL_{a}\times \calL_{b}\times \calL_{c})$-module map
$$\theta:\xymatrix{
\Res_{\fS_{a}\times \fS_{b}\times \fS_{c}}^{\fS_n} (\calL_n) \ar[r] & 
\Res^{\fS_{a+b}}_{\fS_{a}\times\fS_{b}}(\calL_{a+b}) \otimes \calL_{c} \ar[r] &
\calL_{a}\otimes \calL_{b}\otimes \calL_{c},
}$$
and letting $\calL_{a}\otimes \calL_{b}= \largewedge^*\calW$ where
$\calW = 
\left( \calS^{(a-1,1)}\otimes \mathbf{1}_{\fS_{b}}\right)
\oplus 
\left( \mathbf{1}_{\fS_{a}}\otimes \calS^{(b-1,1)}\right)$, 
then 
$$\Res_{\fS_{a}\times \fS_{b}\times \fS_{c}}^{\fS_n} (\calL_n)
= (\largewedge^* \calW\otimes\calL_{c}) 
\oplus \left(\largewedge^* \calW\otimes\calL_{c}\right)^{\oplus2} [1]
\oplus (\largewedge^* \calW\otimes\calL_{c})[2],$$
and
the composite map $\theta:\Res(\calL_n)\rightarrow \calL_{a}\otimes \calL_{b}\otimes\calL_{c}$ is the identity map on $\largewedge^*\calW\otimes \calL_{c}$, and is the zero map on the remaining factors. 

Finally, since $\largewedge^* \calW\otimes\calL_{c} = \calL_{a}\otimes \calL_{b}\otimes\calL_{c} = \calL_{a}\otimes\largewedge^*\calU$, we conclude that $\psi = \theta$, and by Frobenius reciprocity, it follows that 
$d_{\varepsilon_3}\circ d_{\varepsilon_1}=d_{\varepsilon_4}\circ d_{\varepsilon_2}.$
\qed

\subsection{Decategorification via the Frobenius series}
In this section, we show that the chromatic symmetric homology categorifies the chromatic symmetric function.  The decategorification process amounts to computing the bigraded Frobenius series of the homology at $q=t=1$.


\begin{defn}
For $i,j\geq0$, the $(i,j)$th {\em homology} of $G$ is
$$H_{i,j}(G)  = \ker d_{i,j}/\im d_{i+1,j},
\qquad\hbox{so that}\qquad
H_i(G)= \bigoplus_{j\geq0} H_{i,j}(G).$$
The bigraded {\em Frobenius series} of $H_*(G) = \bigoplus_{i,j\geq0}H_{i,j}(G)$ is
$$\mathrm{Frob}_G(q,t) = \sum_{i,j\geq0}(-1)^{i+j} t^i q^j\, \ch \left(H_{i,j}(G)\right),$$
where $\ch: R \rightarrow \Lambda_\bbC: [\calS^\lambda]\mapsto s_\lambda$  is the Frobenius characteristic map from the Grothendieck ring of representations to the ring of symmetric functions defined in equation~\eqref{eqn.frobchar}.
\end{defn}

\begin{lemma}\label{lem.CequalsH} For any graph $G$,
$$\sum_{i,j\geq0} (-1)^{i+j} \ch\left(H_{i,j}(G)\right)
=\sum_{i,j\geq0} (-1)^{i+j}\ch\left(C_{i,j}(G)\right).$$
\end{lemma}
\proof 
This proof is similar to that of the Euler characteristic of chain complexes.  In the following, we supply the proof along with the necessary minor adaptations.  

Any short exact sequence of $\fS_n$-modules $0 \rightarrow A \rightarrow B\rightarrow C\rightarrow 0$ is split exact, so 
$B \cong A\oplus C$ and $\ch(B)= \ch(A)+\ch(C)$.

Consider the chain complex 
$$\xymatrix{
0 \ar[r]
& C_{m,j}(G) \ar[r]^-{d_{m,j}}
& C_{m-1,j}(G) \ar[r] 
& \cdots \ar[r] 
& C_{1,j}(G) \ar[r]^-{d_{1,j}} 
& C_{0,j}(G) \ar[r] 
&0
}$$
for each $j\geq0$.  Let $Z_{i,j}(G) = \ker d_{i,j}$ and $B_{i,j}(G) = \im d_{i+1,j}$.  For $i, j\geq0$, we have short exact sequences
$0 \rightarrow Z_{i,j} \rightarrow C_{i,j} \rightarrow B_{i-1,j} \rightarrow 0$ and
$0 \rightarrow B_{i,j} \rightarrow Z_{i,j} \rightarrow H_{i,j} \rightarrow 0
$,
where $B_{-1,j}$ is understood to be zero.  Thus
\begin{align*}
\ch(C_{i,j}) &= \ch(Z_{i,j}) + \ch(B_{i-1,j})
= \ch(H_{i,j}) + \ch(B_{i,j}) + \ch(B_{i-1,j}).
\end{align*}
Multiplying this by $(-1)^{i+j}$ and summing over all $i,j\geq0$ gives the desired result.
\qed

\begin{theorem}\label{thm.FX} Chromatic symmetric homology categorifies the chromatic symmetric function. 
That is, for any graph $G$, 
$$\Frob_G(1,1) = X_G.$$
\end{theorem}
\proof
Following Lemma~\ref{lem.CequalsH}, we compute
\begin{align*}
\Frob_G(1,1) 
&=\sum_{i,j\geq0} (-1)^{i+j} \ch\left(H_{i,j}(G)\right)\\
&=\sum_{i\geq0} (-1)^i 
	\left(\sum_{j\geq0} (-1)^{j}\ch\left(C_{i,j}(G)\right)\right)
	&\hbox{by Lemma~\ref{lem.CequalsH}}\\
&=\sum_{i\geq0} (-1)^i \sum_{F\subseteq E(G):\atop |F|=i} p_{\lambda(F)} \\  
&= X_G,
		&\hbox{by Equation~\eqref{thm.X_in_p}}.
\end{align*}
\qed

\begin{example}\label{eg.kthree_homology} \label{eg.kthree}
We finish the $K_3$ example by computing its homology and its Frobenius series. Continuing from Example~\ref{eg.kthree_objects},
recall that $\fS_3 = \langle s_1,s_2 \mid s_1^2=e, s_2^2=e, (s_1s_2)^3=e\rangle$.
For convenience, we let $s_\varphi = s_1s_2s_1 = s_2s_1s_2$.
Then $\Ind_{1|2|3} (\calS^1\otimes \calS^1\otimes \calS^1)$ is the regular representation of $\fS_3$, so we choose the (ordered) basis
\begin{align*}
\Ind_{1|2|3} (\calS^1\otimes \calS^1\otimes \calS^1)
& \cong \bbC[\fS_3] \\
& \cong \calS^3\oplus (\calS^{21})^{\oplus2} \oplus \calS^{111}\\
&= \mathrm{span}_\bbC\{e, s_1, s_2, s_2s_1, s_1s_2, s_1s_2s_1\}.
\end{align*}  

We also choose these ordered bases for the following modules:
\begin{align*}
\Ind_{23|1} (\calS^2\otimes \calS^1) 
	&= \bbC[\fS_3](e+s_\varphi)s_1
	= \spn_\bbC\{e+s_2, s_1+s_1s_2, s_2s_1+s_\varphi\},\\
\Ind_{13|2} (\calS^2\otimes \calS^1) 
	&= \bbC[\fS_3](e+s_\varphi)
	= \spn_\bbC\{e+s_\varphi, s_1+s_2s_1, s_2+s_1s_2\},\\
\Ind_{12|3} (\calS^2\otimes \calS^1) 
	&= \bbC[\fS_3](e+s_\varphi)s_2
	=\spn_\bbC\{e+s_1, s_2+s_2s_1, s_1s_2+s_\varphi\}.
\end{align*} 
Also, choose $\calS^3 = \bbC[\fS_3]\sigma = \spn_\bbC\{\sigma\}$, where $\sigma=e+s_1+s_2+s_2s_1+s_1s_2+s_\varphi$ is a Young symmetrizer.
Then every per-edge map in the zeroth grading is an inclusion map:
$$ 
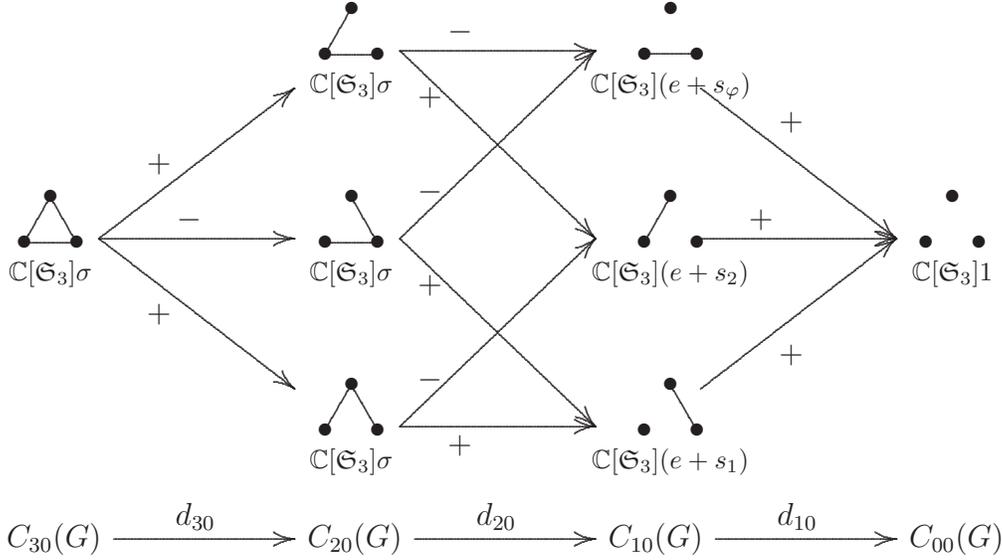
\beginpicture
	\setcoordinatesystem units <4cm,2.5cm>         
	\setplotarea x from 0 to 0, y from -1.9 to 0    

\arrow <8pt> [.2,.67] from 0.1 0 to .75 .8 	\put{$+$} at .3 .4 
\arrow <8pt> [.2,.67] from 0.1 0 to .75 0    \put{$-$} at .4 .1
\arrow <8pt> [.2,.67] from 0.1 0 to .75 -.8  \put{$+$} at .3 -.4

\arrow <8pt> [.2,.67] from 1.1 1 to 1.75 1   \put{$-$}[c] at 1.3 1.1
\arrow <8pt> [.2,.67] from 1.1 1 to 1.75 0   \put{$+$}[c] at 1.2 .75
\arrow <8pt> [.2,.67] from 1.1 0 to 1.75 1   \put{$-$}[c] at 1.2 .25
\arrow <8pt> [.2,.67] from 1.1 0 to 1.75 -1  \put{$+$}[c] at 1.2 -.25
\arrow <8pt> [.2,.67] from 1.1 -1 to 1.75 0  \put{$-$}[c] at 1.2 -.75
\arrow <8pt> [.2,.67] from 1.1 -1 to 1.75 -1 \put{$+$}[c] at 1.3 -1.1

\arrow <8pt> [.2,.67] from 2.1 .8 to 2.75 0  \put{$+$}[c] at 2.4 .62
\arrow <8pt> [.2,.67] from 2.1 0 to 2.75 0   \put{$+$}[c] at 2.3 .1
\arrow <8pt> [.2,.67] from 2.1 -.8 to 2.75 0 \put{$+$}[c] at 2.4 -.62 

\put{$\beginpicture
	\setcoordinatesystem units <.7cm,.7cm>         
	\setplotarea x from 0 to 2, y from 0 to 0    
	\multiput{$\bullet$} at 0 0 1 0 0.5 0.867 /
	\plot 0 0 1 0 0.5 .867 0 0 /
	\put{\footnotesize$\bbC[\fS_3]\sigma$} at .5 -.6
	\endpicture$}[c] at 0 0

\put{$\beginpicture
	\setcoordinatesystem units <.7cm,.7cm>         
	\setplotarea x from 0 to 2, y from 0 to 0    
	\multiput{$\bullet$} at 0 0 1 0 0.5 0.867 /
	\plot 1 0 0 0 0.5 .867 /
	\put{\footnotesize$\bbC[\fS_3]\sigma$} at .5 -.6
	\endpicture$}[c] at 1 1
\put{$\beginpicture
	\setcoordinatesystem units <.7cm,.7cm>         
	\setplotarea x from 0 to 2, y from 0 to 0    
	\multiput{$\bullet$} at 0 0 1 0 0.5 0.867 /
	\plot 0.5 .867 1 0 0 0 /
	\put{\footnotesize$\bbC[\fS_3]\sigma$} at .5 -.6
	\endpicture$}[c] at 1 0	
\put{$\beginpicture
	\setcoordinatesystem units <.7cm,.7cm>         
	\setplotarea x from 0 to 2, y from 0 to 0    
	\multiput{$\bullet$} at 0 0 1 0 0.5 0.867 /
	\plot 0 0 0.5 .867 1 0 /
	\put{\footnotesize$\bbC[\fS_3]\sigma$} at .5 -.6
	\endpicture$}[c] at 1 -1

\put{$\beginpicture
	\setcoordinatesystem units <.7cm,.7cm>         
	\setplotarea x from 0 to 2, y from 0 to 0    
	\multiput{$\bullet$} at 0 0 1 0 0.5 0.867 /
	\plot 0 0 1 0 /
	\put{\footnotesize$\bbC[\fS_3](e+s_\varphi)$} 
		at .5 -.6
	\endpicture$}[c] at 2 1
\put{$\beginpicture
	\setcoordinatesystem units <.7cm,.7cm>         
	\setplotarea x from 0 to 2, y from 0 to 0    
	\multiput{$\bullet$} at 0 0 1 0 0.5 0.867 /
	\plot 0 0 0.5 .867 /
	\put{\footnotesize$\bbC[\fS_3](e+s_2)$} 
		at .5 -.6
	\endpicture$}[c] at 2 0
\put{$\beginpicture
	\setcoordinatesystem units <.7cm,.7cm>         
	\setplotarea x from 0 to 2, y from 0 to 0    
	\multiput{$\bullet$} at 0 0 1 0 0.5 0.867 /
	\plot 1 0 0.5 .867 /
	\put{\footnotesize$\bbC[\fS_3](e+s_1)$} 
		at .5 -.6
	\endpicture$}[c] at 2 -1

\put{$\beginpicture
	\setcoordinatesystem units <.7cm,.7cm>         
	\setplotarea x from 0 to 2, y from 0 to 0    
	\multiput{$\bullet$} at 0 0 1 0 0.5 0.867 /
	\put{\footnotesize$\bbC[\fS_3]1$} at .5 -.6
	\endpicture$} at 3 0
\put{$C_{30}(G)$}[c] at -.05 -1.6
\put{$C_{20}(G)$}[c] at .95 -1.6
\put{$C_{10}(G)$}[c] at 1.95 -1.6
\put{$C_{00}(G)$}[c] at 2.95 -1.6
\arrow <8pt> [.2,.67] from 0.15 -1.6 to 0.75 -1.6 
	\put{$d_{30}$} at .42 -1.48
\arrow <8pt> [.2,.67] from 1.15 -1.6 to 1.75 -1.6 	
	\put{$d_{20}$} at 1.42 -1.48
\arrow <8pt> [.2,.67] from 2.15 -1.6 to 2.75 -1.6 	
	\put{$d_{10}$} at 2.42 -1.48
\endpicture
$$
Thus $\ker d_{30} \cong 0$, $\im d_{30} = \ker d_{20} \cong \calS^3$, 
$\im d_{20} \cong  (\calS^3)^{\oplus2}$,  $\ker d_{10} \cong (\calS^3)^{\oplus2}\oplus \calS^{21}$,
$\im d_{10} \cong \calS^3 \oplus (\calS^{21})^{\oplus2}$,
and $\ker d_{00} \cong \bbC[\fS_3]$.

In the first grading, we choose the ordered basis 
$\{e-s_\varphi, s_1-s_2s_1, s_2-s_1s_2 \}$ for $\Ind_{13|2}(\calS^{11}\otimes \calS^1) = \bbC[\fS_3](e-s_\varphi)$.  
It follows that 
$\Ind_{23|1} (\calS^{11}\otimes \calS^1) = \bbC[\fS_3](e-s_\varphi)s_1$
and $\Ind_{12|3} (\calS^{11}\otimes \calS^1) = \bbC[\fS_3](e-s_\varphi)s_2$.
Let $\tau = (s_2+s_2s_1)(e-s_\varphi)$ be a Young symmetrizer in $\Ind_{13|2} (\calS^{11}\otimes \calS^2)$, and choose $\calS^{21}=\bbC[\fS_3]\tau$.
Then the per-edge maps $d_{*11}, d_{1*1}, d_{11*}, d_{0*1}, d_{*01}$ are inclusion maps (up to sign).  Also, $d_{01*}(\tau) = d_{*10}(\tau) =  \tau s_1$, while
$d_{10*}(\tau) =d_{1*0}(\tau) =  \tau s_2$.  
$$ 
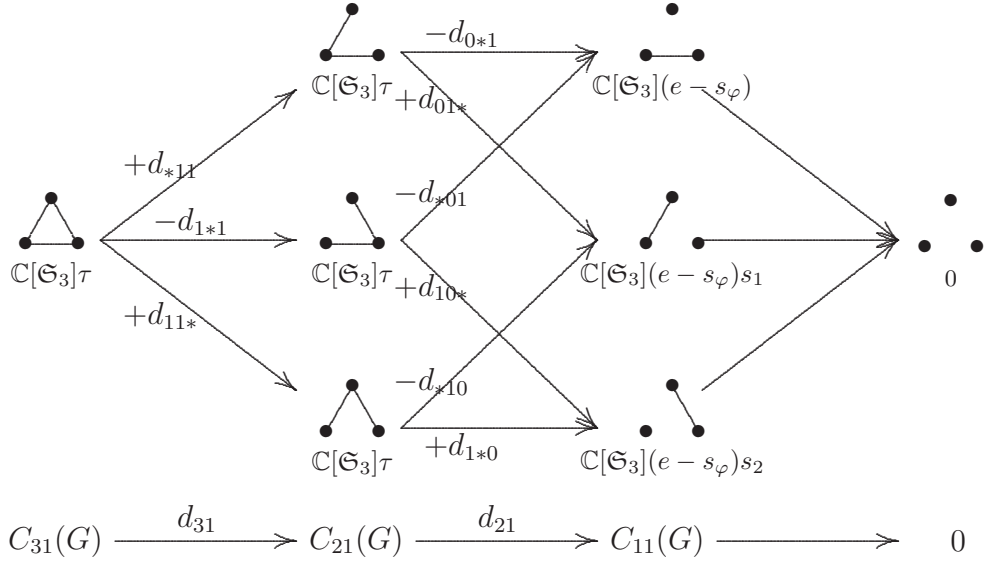
\beginpicture
	\setcoordinatesystem units <4cm,2.5cm>         
	\setplotarea x from 0 to 0, y from -2 to 1.1    

\arrow <8pt> [.2,.67] from 0.1 0 to .75 .8 	\put{$+d_{*11}$} at .3 .4 
\arrow <8pt> [.2,.67] from 0.1 0 to .75 0    \put{$-d_{1*1}$} at .4 .1
\arrow <8pt> [.2,.67] from 0.1 0 to .75 -.8  \put{$+d_{11*}$} at .3 -.4

\arrow <8pt> [.2,.67] from 1.1 1 to 1.75 1   \put{$-d_{0*1}$}[c] at 1.3 1.1
\arrow <8pt> [.2,.67] from 1.1 1 to 1.75 0   \put{$+d_{01*}$}[c] at 1.2 .75
\arrow <8pt> [.2,.67] from 1.1 0 to 1.75 1   \put{$-d_{*01}$}[c] at 1.2 .25
\arrow <8pt> [.2,.67] from 1.1 0 to 1.75 -1  \put{$+d_{10*}$}[c] at 1.2 -.25
\arrow <8pt> [.2,.67] from 1.1 -1 to 1.75 0  \put{$-d_{*10}$}[c] at 1.2 -.75
\arrow <8pt> [.2,.67] from 1.1 -1 to 1.75 -1 \put{$+d_{1*0}$}[c] at 1.3 -1.1

\arrow <8pt> [.2,.67] from 2.1 .8 to 2.75 0  
\arrow <8pt> [.2,.67] from 2.1 0 to 2.75 0   
\arrow <8pt> [.2,.67] from 2.1 -.8 to 2.75 0 
 
\put{$\beginpicture
	\setcoordinatesystem units <.7cm,.7cm>         
	\setplotarea x from 0 to 2, y from 0 to 0    
	\multiput{$\bullet$} at 0 0 1 0 0.5 0.867 /
	\plot 0 0 1 0 0.5 .867 0 0 /
	\put{\footnotesize$\bbC[\fS_3]\tau$} at .5 -.6
	\endpicture$}[c] at 0 0

\put{$\beginpicture
	\setcoordinatesystem units <.7cm,.7cm>         
	\setplotarea x from 0 to 2, y from 0 to 0    
	\multiput{$\bullet$} at 0 0 1 0 0.5 0.867 /
	\plot 1 0 0 0 0.5 .867 /
	\put{\footnotesize$\bbC[\fS_3]\tau$} at .5 -.6
	\endpicture$}[c] at 1 1
\put{$\beginpicture
	\setcoordinatesystem units <.7cm,.7cm>         
	\setplotarea x from 0 to 2, y from 0 to 0    
	\multiput{$\bullet$} at 0 0 1 0 0.5 0.867 /
	\plot 0.5 .867 1 0 0 0 /
	\put{\footnotesize$\bbC[\fS_3]\tau$} at .5 -.6
	\endpicture$}[c] at 1 0	
\put{$\beginpicture
	\setcoordinatesystem units <.7cm,.7cm>         
	\setplotarea x from 0 to 2, y from 0 to 0    
	\multiput{$\bullet$} at 0 0 1 0 0.5 0.867 /
	\plot 0 0 0.5 .867 1 0 /
	\put{\footnotesize$\bbC[\fS_3]\tau$} at .5 -.6
	\endpicture$}[c] at 1 -1

\put{$\beginpicture
	\setcoordinatesystem units <.7cm,.7cm>         
	\setplotarea x from 0 to 2, y from 0 to 0    
	\multiput{$\bullet$} at 0 0 1 0 0.5 0.867 /
	\plot 0 0 1 0 /
	\put{\footnotesize$\bbC[\fS_3](e-s_\varphi)$} 
		at .5 -.6
	\endpicture$}[c] at 2 1
\put{$\beginpicture
	\setcoordinatesystem units <.7cm,.7cm>         
	\setplotarea x from 0 to 2, y from 0 to 0    
	\multiput{$\bullet$} at 0 0 1 0 0.5 0.867 /
	\plot 0 0 0.5 .867 /
	\put{\footnotesize$\bbC[\fS_3](e-s_\varphi)s_1$} 
		at .5 -.6
	\endpicture$}[c] at 2 0
\put{$\beginpicture
	\setcoordinatesystem units <.7cm,.7cm>         
	\setplotarea x from 0 to 2, y from 0 to 0    
	\multiput{$\bullet$} at 0 0 1 0 0.5 0.867 /
	\plot 1 0 0.5 .867 /
	\put{\footnotesize$\bbC[\fS_3](e-s_\varphi)s_2$} 
		at .5 -.6
	\endpicture$}[c] at 2 -1

\put{$\beginpicture
	\setcoordinatesystem units <.7cm,.7cm>         
	\setplotarea x from 0 to 2, y from 0 to 0    
	\multiput{$\bullet$} at 0 0 1 0 0.5 0.867 /
	\put{\scriptsize$0$} at .5 -.6
	\endpicture$} at 3 0
\put{$C_{31}(G)$}[c] at -.05 -1.6
\put{$C_{21}(G)$}[c] at .95 -1.6
\put{$C_{11}(G)$}[c] at 1.95 -1.6
\put{$0$}[c] at 2.95 -1.6
\arrow <8pt> [.2,.67] from 0.15 -1.6 to 0.75 -1.6 
	\put{$d_{31}$} at .42 -1.48
\arrow <8pt> [.2,.67] from 1.15 -1.6 to 1.75 -1.6 	
	\put{$d_{21}$} at 1.42 -1.48
\arrow <8pt> [.2,.67] from 2.15 -1.6 to 2.75 -1.6 	
\endpicture
$$
Thus $\ker d_{31} \cong 0$, $\im d_{31} = \ker d_{21} \cong \calS^{21}$, 
$\im d_{21} = (\calS^{21})^{\oplus2}$,
and $\ker d_{11}\cong  (\calS^{21}\oplus\calS^{111})^{\oplus3}$.

Lastly, in the second grading, let $\upsilon = e-s_1-s_2 + s_2s_1+s_1s_2- s_\varphi$, so that $\calS^{111}=\bbC[\fS_3]\upsilon =\spn_\bbC\{\upsilon\}$.  Then
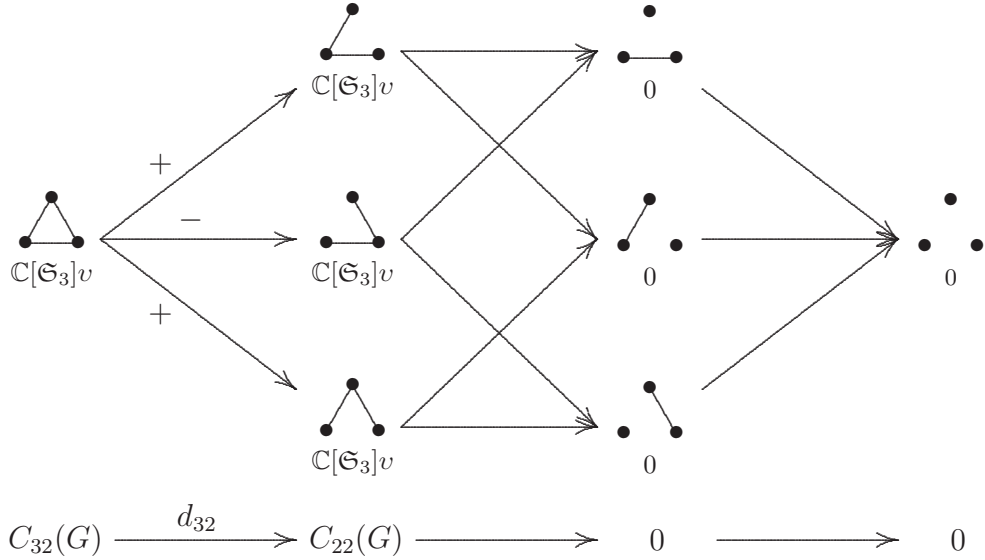
$$ 
\beginpicture
	\setcoordinatesystem units <4cm,2.5cm>         
	\setplotarea x from 0 to 0, y from -2 to 1.1    

\arrow <8pt> [.2,.67] from 0.1 0 to .75 .8 	\put{$+$} at .3 .4 
\arrow <8pt> [.2,.67] from 0.1 0 to .75 0    \put{$-$} at .4 .1
\arrow <8pt> [.2,.67] from 0.1 0 to .75 -.8  \put{$+$} at .3 -.4

\arrow <8pt> [.2,.67] from 1.1 1 to 1.75 1   
\arrow <8pt> [.2,.67] from 1.1 1 to 1.75 0   
\arrow <8pt> [.2,.67] from 1.1 0 to 1.75 1   
\arrow <8pt> [.2,.67] from 1.1 0 to 1.75 -1  
\arrow <8pt> [.2,.67] from 1.1 -1 to 1.75 0  
\arrow <8pt> [.2,.67] from 1.1 -1 to 1.75 -1 

\arrow <8pt> [.2,.67] from 2.1 .8 to 2.75 0  
\arrow <8pt> [.2,.67] from 2.1 0 to 2.75 0   
\arrow <8pt> [.2,.67] from 2.1 -.8 to 2.75 0 

\put{$\beginpicture
	\setcoordinatesystem units <.7cm,.7cm>         
	\setplotarea x from 0 to 2, y from 0 to 0    
	\multiput{$\bullet$} at 0 0 1 0 0.5 0.867 /
	\plot 0 0 1 0 0.5 .867 0 0 /
	\put{\footnotesize$\bbC[\fS_3]\upsilon$} at .5 -.6
	\endpicture$}[c] at 0 0

\put{$\beginpicture
	\setcoordinatesystem units <.7cm,.7cm>         
	\setplotarea x from 0 to 2, y from 0 to 0    
	\multiput{$\bullet$} at 0 0 1 0 0.5 0.867 /
	\plot 1 0 0 0 0.5 .867 /
	\put{\footnotesize$\bbC[\fS_3]\upsilon$} at .5 -.6
	\endpicture$}[c] at 1 1
\put{$\beginpicture
	\setcoordinatesystem units <.7cm,.7cm>         
	\setplotarea x from 0 to 2, y from 0 to 0    
	\multiput{$\bullet$} at 0 0 1 0 0.5 0.867 /
	\plot 0.5 .867 1 0 0 0 /
	\put{\footnotesize$\bbC[\fS_3]\upsilon$} at .5 -.6
	\endpicture$}[c] at 1 0	
\put{$\beginpicture
	\setcoordinatesystem units <.7cm,.7cm>         
	\setplotarea x from 0 to 2, y from 0 to 0    
	\multiput{$\bullet$} at 0 0 1 0 0.5 0.867 /
	\plot 0 0 0.5 .867 1 0 /
	\put{\footnotesize$\bbC[\fS_3]\upsilon$} at .5 -.6
	\endpicture$}[c] at 1 -1

\put{$\beginpicture
	\setcoordinatesystem units <.7cm,.7cm>         
	\setplotarea x from 0 to 2, y from 0 to 0    
	\multiput{$\bullet$} at 0 0 1 0 0.5 0.867 /
	\plot 0 0 1 0 /
	\put{\footnotesize$0$} 
		at .5 -.6
	\endpicture$}[c] at 2 1	
\put{$\beginpicture
	\setcoordinatesystem units <.7cm,.7cm>         
	\setplotarea x from 0 to 2, y from 0 to 0    
	\multiput{$\bullet$} at 0 0 1 0 0.5 0.867 /
	\plot 0 0 0.5 .867 /
	\put{\footnotesize$0$} 
		at .5 -.6
	\endpicture$}[c] at 2 0
\put{$\beginpicture
	\setcoordinatesystem units <.7cm,.7cm>         
	\setplotarea x from 0 to 2, y from 0 to 0    
	\multiput{$\bullet$} at 0 0 1 0 0.5 0.867 /
	\plot 1 0 0.5 .867 /
	\put{\footnotesize$0$} 
		at .5 -.6
	\endpicture$}[c] at 2 -1

\put{$\beginpicture
	\setcoordinatesystem units <.7cm,.7cm>         
	\setplotarea x from 0 to 2, y from 0 to 0    
	\multiput{$\bullet$} at 0 0 1 0 0.5 0.867 /
	\put{\scriptsize$0$} at .5 -.6
	\endpicture$} at 3 0

\put{$C_{32}(G)$}[c] at -.05 -1.6
\put{$C_{22}(G)$}[c] at .95 -1.6
\put{$0$}[c] at 1.95 -1.6
\put{$0$}[c] at 2.95 -1.6
\arrow <8pt> [.2,.67] from 0.15 -1.6 to 0.75 -1.6
	\put{$d_{32}$} at .42 -1.48
\arrow <8pt> [.2,.67] from 1.15 -1.6 to 1.75 -1.6 	
\arrow <8pt> [.2,.67] from 2.15 -1.6 to 2.75 -1.6 	
\endpicture
$$
So $\ker d_{32}\cong 0$, $\im d_{32} \cong \calS^{111}$, and $\ker d_{32}=(\calS^{111})^{\oplus2}$.

In summary, the bigraded homology of $G=K_3$ is
$$\begin{array}{c|ccc}
2 
	& (\calS^{111})^{\oplus2} & &  \\
1 
	& 0
	& \calS^{21}\oplus(\calS^{111})^{\oplus3} & \\
0 
 	& 0 
	& \calS^{21}
	& \calS^{111} \\
\hline 
& H_2(K_3)& H_1(K_3)& H_0(K_3)
\end{array}$$
and the bigraded Frobenius series is
\begin{align*}
\Frob_{K_3}(q,t) 
&= 2q^2t^2s_{111} +qt\left(s_{21}+3s_{111} \right)-ts_{21}+s_{111}\\
&= (1-q)ts_{21} +(1+qt)(1+2qt)s_{111}\\
&= (1-q)tm_{21}+ (1-2t+5qt+2q^2t^2)m_{111}.
\end{align*}
The specializations are $\Frob_{K_3}(1,t) = (1+t)(1+2t)s_{111}$,
and $\Frob_{K_3}(1,1) = 6s_{111}$.
\end{example}

\section{Properties of $H_*(G)$}\label{sec.properties}

\subsection{Loops and multiple edges}
\begin{prop}\label{prop.loop}  If $G$ contains a loop, then $H_*(G)=0$.
\end{prop}
\proof 
Suppose $G$ has $m$ edges $\{e_1,\ldots,e_m\}$, and $e_m=\ell$ is a loop in $G$.  
Let $(U_*, d^U_*)$ be the chain complex obtained by restricting $(C_*(G), d_*)$ to the modules indexed by states which contain $\ell$, and let $A_* = U_*[1]$ be the first translate of $U_*$ so that $A_i = U_{i+1}$ for all $i\geq0$.  Let $(B_*, d^B_*)$ be the chain complex obtained by restricting $(C_*(G), d_*)$ to the modules indexed by states which do not contain $\ell$.

Observe that the only per-edge maps from modules in $A_*$ to modules in $B_*$ are of the form $\calM_F \rightarrow \calM_{F-\ell}$, and each of these maps is the identity map as $\calM_F  = \calM_{F-\ell}$ if $\ell$ is a loop.  Moreover, $A_*=B_*$ and $d^A_*= -d^B_*$ (by the definition of the signed per-edge maps), so we have
$$\xymatrix{
\ar[r]& 
	A_{i,j} \ar[r]^{-d^{B}_{i,j}} \ar[rd]^{\mathrm{id}}& 
	A_{i-1,j} \ar[r]^{-d^{B}_{i-1,j}} \ar[rd]^{\mathrm{id}}&
	A_{i-2,j} \ar[r]& \\
\ar[r]&
	B_{i+1,j} \ar[r]_{d^{B}_{i+1,j}} &
	B_{i,j} \ar[r]_{d^{B}_{i,j}}&
	B_{i-1,j} \ar[r]&
}$$
where
$$C_{i,j}(G) = A_{i-1,j}\oplus B_{i,j} \qquad\hbox{and} \qquad
d_{i,j}= \begin{bmatrix}
-d_{i-1,j}^B &0\\ 
\mathrm{id} &d_{i,j}^B
\end{bmatrix}.$$
In other words, $C_*(G)$ is the mapping cone of the negative identity map 
$-\mathrm{id}:A_*\rightarrow B_*$, so $C_*(G)$ is exact by Exercise 1.5.1 in~\cite{weibel1995}.

We supply the details below.  Since $d_{i,j}(a,b) = \left(-d_{i-1,j}^B(a), a+d_{i,j}^B(b)\right)$, then
$$\ker d_{i,j} = \left\{ (a,b) \mid a=-d_{i,j}^B(b)\in \im d_{i,j}^B \subseteq \ker d_{i-1,j}^B \right\}
=\left\{ \left(-d_{i,j}^B(b), b\right)\mid b\in B_{i,j}\right\}.
$$
On the other hand,
$\im d_{i+1} = \left\{ (-d_i^B(x),x+d_{i+1}^B(y)) \mid (x,y)\in A_i\oplus B_{i+1}\right\}$.  Since $\im d_{i+1}^B \subseteq \ker d_i^B$, then $(0,d_{i+1}^B(y)) = (-d_i^B(z),z)$ for some $z\in \ker d_i^B \subseteq B_i$.
Thus 
$$\im d_{i+1} = \left\{ \left(-d_i^B(x+z), x+z\right) \mid x+z\in B_i \right\} = \ker d_i.$$
Therefore, $C_*(G)$ is exact, and $H_*(G)=0$. \qed

\begin{prop}  Let $G$ be a multigraph with edges $e$ and $e'$ each incident to the vertices $x$ and $y$. Then $H_*(G) = H_*(G-e')$.
\end{prop}
\proof Suppose $G$ has $m$ edges $\{e_1,\ldots, e_m\}$ and $e_m=e$.
Let $Z=G-e'$.  Since $Z$ is a subgraph of $G$, then $C_{i}(Z)\subseteq C_{i}(G)$ for all $i\geq0$.  Define $C_{i}(G,Z)= C_{i}(G)/C_{i}(Z)$.  As the differential $d_{i}: C_{i}(G) \rightarrow C_{i-1}(G)$ sends $C_{i}(Z)$ to $C_{i-1}(Z)$,
then there is an induced complex
$$\xymatrix{
\cdots \ar[r]&
C_{i}(G,Z)\ar[r]^-{d_{i}}&
C_{i-1}(G,Z)\ar[r]&\cdots.
} $$
Define $H_{i}(G,Z) = \ker d_{i}/\im d_{i+1}$. 
There are short exact sequences
$$\xymatrix{0\ar[r]&
C_i(Z)\ar@{^{(}->}[r]& 
C_i(G)\ar@{>>}[r]&
C_i(G,Z)\ar[r]&0}$$
for each $i\geq0$, and these induce
the long exact sequence of homology
$$\xymatrix{
\cdots\ar[r]&
H_i(Z)\ar[r]&
H_i(G)\ar[r]&
H_i(G,Z)\ar[r]&
H_{i-1}(Z)\ar[r]&\cdots .
}$$
We shall show that $H_{i,j}(G,Z)=0$ for all $i,j\geq0$, from whence it follows that $H_*(G) = H_*(Z)$.

To compute relative homology $H_*(G,Z)$,  we only need to consider the states in $Q(G)$ which contain $e'$. Let $(A_*, d_*^A)$ denote the chain complex obtained by restricting $(C_*(G), d_*)$ to the modules indexed by the states which contain both $e'$ and $e$, and let $(B_*, d_*^B)$ denote the chain complex obtained by restricting $(C_*(G), d_*)$ to the modules indexed by the states which contain $e'$ but not $e$.  The only per-edge maps from modules in $A_*$ to modules in $B_*$ are of the form $\calM_F \rightarrow \calM_{F-e}$, and each of these maps is the identity map because $e'\in F$ implies $\calM_F = \calM_{F-e}$.
Similar to the case of loops, we see from this that $C_*(G,Z)$ is the mapping cone of the negative identity map $-\mathrm{id}:A_*[1] \rightarrow B_*$, so $C_*(G,Z)$ is exact, and hence $H_*(G,Z)=0$.\qed

\subsection{Disjoint unions}
Given graphs $A$ and $B$, let $A+B$ denote their disjoint union.  Since the power sum symmetric functions are multiplicative, one can deduce from this that $X_{A+B} = X_A X_B$. The lift to homology is the next result.  The specialization of the Frobenius series at $q=t=1$ recovers the multiplicative property of the chromatic symmetric function of a disjoint union of graphs.

\begin{prop}\label{prop.disjointunion}
For $i,j\geq0$, 
$$H_{i,j}(A+B) = \bigoplus_{p+r=i \atop q+s=j} \Ind_{\fS_A\times \fS_B}^{\fS_{A\cup B}} \left( H_{p,q}(A)\otimes H_{r,s}(B)\right).$$
\end{prop}
\proof

This is essentially given by the K\"unneth formula.
Let $G = A+B$, and suppose we enumerate the edges of $G$ by enumerating the edges of $B$ first, then followed by the edges of $A$.  Now, $Q(G) \cong Q(A)\times Q(B)$, so for each $i\geq0$,
$$C_i(G) = \bigoplus_{F\in Q(G):\atop |F|=i} \calM_F 
= \bigoplus_{(F_A, F_B)\in Q(A)\times Q(B):\atop
|F_A|+|F_B|=i }\Ind_{\fS_A\times \fS_B}^{\fS_{A\cup B}} \left(\calM_{F_A}\otimes \calM_{F_B} \right).$$
For the moment, we consider the $(\fS_A\times \fS_B)$-modules
$$\widetilde{C}_i(G)
= \bigoplus_{(F_A, F_B)\in Q(A)\times Q(B):\atop
|F_A|+|F_B|=i} \left(\calM_{F_A}\otimes \calM_{F_B} \right),$$
so that
$\widetilde{C}_{i,j}(G)
= \bigoplus_{p+r=i\atop q+s=j} 
\left(C_{p,q}(A)\otimes C_{r,s}(B) \right)$, and the chain maps are given by
$$\widetilde{d}_{i,j}(c_A \otimes c_B) = d_{p,q}^A(c_A) \otimes c_B + (-1)^{p} c_A \otimes d_{r,s}^B(c_B).$$
These modules are free and hence are flat, so by~\cite[Theorem 3.6.3]{weibel1995}, there is an exact sequence
$$ 0\rightarrow 
\bigoplus_{p+r=i\atop q+s=j} \left(\widetilde{H}_{p,q}(A) \otimes \widetilde{H}_{r,s}(B)\right) \rightarrow \widetilde{H}_{i,j}(G)
\rightarrow 
\bigoplus_{p+r=i-1\atop q+s=j} \mathrm{Tor}_1\left(\widetilde{H}_{p,q}(A), \widetilde{H}_{r,s}(B)\right) 
\rightarrow0.$$
Moreover, there is zero torsion, again due to freeness.  Therefore,
\begin{equation}\label{eqn.tensor} 
\widetilde{H}_{i,j}(G) \cong \bigoplus_{p+r=i\atop q+s=j}  
\left(\widetilde{H}_{p,q}(A) \otimes \widetilde{H}_{r,s}(B)\right). 
\end{equation}
Lastly, induce these to $\fS_{A\cup B}$-modules to get the desired result. \qed

\begin{corollary}
$\Frob_{A+B}(q,t) = \Frob_A (q,t)\cdot \Frob_B(q,t)$.
\qed
\end{corollary}

\begin{example} If $G$ is a graph with $n$ vertices, then the homology of the disjoint union of $G$ with a single vertex is obtained by branching.  
If $H_{i,j}(G) = \bigoplus_{\lambda} \left(\calS^{\lambda}\right)^{\oplus m_\lambda}$, then 
$$H_{i,j}(G+\bullet) = \Ind_{\fS_n}^{\fS_{n+1}} H_{i,j}(G)
\cong \bigoplus_{\mu} \left(\calS^{\mu}\right)^{\oplus m_\lambda},$$
where the sum is over all partitions $\mu$ which can be obtained by adding a box to the partitions $\lambda$ indexing the irreducible factors of $H_{i,j}(G)$.
\end{example}

\section{Categorifying the recursive formula for $X_G$}\label{sec.mv}

If the graph $G$ has $n$ vertices, then $X_G$ is homogeneous of degree $n$, so it follows that $X_G$ does not satisfy a deletion-contraction type recurrence.  However, if $G$ contains triangles, then the following result provides a method for expressing $X_G$ as a linear combination of chromatic symmetric functions of graphs with fewer edges.
\begin{theorem}[{\cite[Proposition 3.1]{guaypaquet2013} and \cite[Theorem 3.1]{os2014}}]
\label{thm.dd}
Let $G$ be a graph where the edges $e_1,e_2,e_3\in E(G)$ form a triangle.  Then
$$X_G  = X_{G-e_1} + X_{G-e_2} - X_{G-e_1-e_2}.$$
\qed
\end{theorem}

The homological analogue of this three term recurrence equation is 
a Mayer-Vietoris type long exact sequence.
Let $A=G-e_1$, $B=G-e_2$, and let
\begin{equation}\label{eqn.CAUB}
C_i(A\cup B) =  \bigoplus_{|F|=i, \atop F\in Q(A) \cup Q(B)} \calM_F
\end{equation}
be the submodule of $C_i(G)$ that is the direct sum of the modules $\calM_F$ where $F$ is a state of the subgraph $A$ or a state of the subgraph $B$.
The differential $d:C_{i}(G) \rightarrow C_{i-1}(G)$ sends $C_i(A\cup B)$ to $C_{i-1}(A\cup B)$, so $C_*(A\cup B)$ forms a chain complex.  Let $H_*(A\cup B)$ be the bigraded homology resulting from $C_*(A\cup B)$.

\begin{lemma} \label{lem.modchain}
$H_*(G) \cong H_*(A\cup B)$.
\end{lemma}

\proof
We shall show that the inclusion $\iota: C_*(A\cup B) \hookrightarrow C_*(G)$ induces isomorphisms on homology.  That is, we shall define a chain map $f_*:C_*(G) \rightarrow C_*(A\cup B)$ and prism operators $s_i:C_i(G) \rightarrow C_{i+1}(G)$
$$\xymatrix{
\cdots\ar[r]& C_{i+1}(G) \ar[r]^{d_{i+1}} \ar[d]_{\iota\circ f} & 
		C_{i}(G) \ar[r]^{d_i} \ar[d]_{\iota\circ f} \ar[ld]_{s_i} & 
		C_{i-1}(G) \ar[r] \ar[d]_{\iota\circ f} \ar[ld]_{s_{i-1}} &
		\cdots\\
\cdots\ar[r]& C_{i+1}(G) \ar[r]^{d_{i+1}} & 
		C_{i}(G) \ar[r]^{d_i} & 
		C_{i-1}(G) \ar[r] & \cdots\\
}$$
such that
$\mathrm{id}_{C_*(G)} - \iota\circ f = ds+sd$ and $f\circ \iota = \mathrm{id}_{C_*(A\cup B)}$, whence it follows that $H_*(G)\cong H_*(A\cup B)$.

To define $f_*$, partition the states of $G$ in the following way.
Suppose $G$ has $m$ edges $\{e_1,\ldots, e_m\}$ and $e_1,e_2,e_3$ form a triangle $T$.  Then $G$ has $2^{m-3}$ states containing $T$.  Each such state $X$ indexes a Boolean sublattice of $Q(G)$ isomorphic to $2^T$: 
\begin{equation}\label{eqn.X}
\xymatrix{
 & X_1 \ar[r] \ar[dr] 
 & X_{12} \ar[dr] & \\
X \ar[ur] \ar[r] \ar[dr] 
	& X_2 \ar[ur] \ar[dr] 
	& X_{13} \ar[r]
	& X_{123} \\
 & X_3 \ar[ur] \ar[r] 
 & X_{23} \ar[ur]& \\
}
\end{equation}
where $X_i = X-e_i$, $X_{ij} = X-e_i-e_j$, and $X_{123} = X-T$. We note that $Q(G) \cong 2^T \times 2^{m-3}$ as Boolean lattices, and this observation will be important for the computations to come.

As the differentials $d_*$ are defined on the modules indexed by the states of $G$, we can do the same to define $f_*$ and $s_*$.  Partition the states of $G$ into $2^{m-3}$ sublattices, each indexed by an $X\supseteq T$.
On each such sublattice, define $f_*$ on the eight corresponding modules as follows.
\begin{enumerate}
\item[(i)] If $F = X$,  let $f=0$ on $\calM_F$.
\item[(ii)] The restriction of $f$ to $\calM=\calM_{X_1}\oplus \calM_{X_2}\oplus \calM_{X_3}$ is given by
$$f:\calM\rightarrow \calM_{X_1}\oplus\calM_{X_2} : (x,y,z)\mapsto (x-z,y+z). $$
\item[(iii)] If $F$ is one of $X_{12}, X_{13}, X_{23}$ or $X_{123}$, let $f$ be the identity map on $\calM_F$.
\end{enumerate}
We remark that since $T$ is a triangle, then deleting any one of the edges  $e_1, e_2$ or $e_3$ will not disconnect $X$, and therefore, $\calM_X = \calM_{X_i}$ for $i=1,2,3$.

We now check that $f_*$ is a chain map. 
Suppose $X=\{e_1,e_2,e_3,e_{i_1}, \ldots, e_{i_r}\}$, and consider the restriction of $f$ to $\calM_X$. Then
$$\xymatrix{
\calM_X\ar[d]_{f_{r+3}=0} \ar[r]^-{d_{r+3}}& 
	\calM_{X_1}\oplus\calM_{X_2} \oplus \calM_{X_3}	\oplus \calN 
		\ar[d]^{f_{r+2}} \\
0\ar[r]^-{d_{r+3}}& 
	\calM_{X_1}\oplus\calM_{X_2}  \\
}$$
where $\calN=\bigoplus_{j=1}^r \calM_{X_{i_j}}$.  For $1\leq j\leq r$, each $X_{i_j}\supseteq T$, so $f|_{X_{i_j}}=0$.
Depending on the parity of $m$, there are two different possibilities of signs. 
So for any $x\in\calM_X$,
$$f_{r+2}(d_{r+3}(x)) = f_{r+2}(\pm x,\mp x, \pm x,n)=0 = d_{r+3}(f_{r+3}(x)),$$
where $n\in\calN$ is the image of $x$ under $d_{r+3}$.

Next, consider the restriction of $f$ to $\calM_{X_3}$.  Then
$$\xymatrix{
\calM_{X_3}
	\ar[r]^-{d_{r+2}}\ar[d]_{f_{r+2}} &
\calM_{X_{13}}\oplus\calM_{X_{23}} \oplus \calN
	\ar[d]^{f_{r+1}} \\  
\calM_{X_1}\oplus \calM_{X_2}
	\ar[r]^-{d_{r+2}}&
\calM_{X_{12}}\oplus\calM_{X_{13}}\oplus\calM_{X_{23}}\oplus \calP\\
}$$
where $\calN=\bigoplus_{j=1}^r \calM_{X_{3,i_j}}$ is a direct sum of modules such that each summand is indexed by a state $F'$ that contains $e_1$ and $e_2$ but not $e_3$, and $\calP = \bigoplus_{j=1}^r (\calM_{X_{1,i_j}}\oplus \calM_{X_{2,i_j}})$.
Observe that $\calM_{X_{3,i_j}}=\calM_{X_{1,i_j}}=\calM_{X_{2,i_j}}$ since they have the same connected graph components, so that when $f$ is restricted to $\calM_{X_{3,i_j}}$, then by definition, 
$$f:\calM_{X_{3,i_j}} \rightarrow \calM_{X_{1,i_j}}\oplus \calM_{X_{2,i_j}}:z \mapsto (-z,z).$$
Therefore, for any $x\in \calM_{X_3}$,
\begin{align*}
d_{r+2}(f_{r+2}(x))
&= d_{r+2}(-x, x)\\
&= \left((\mp(-x))^{\mu_1} +(\mp x)^{\mu_1}, 
	(\pm(-x))^{\mu_2}, 
	(\pm x)^{\mu_3},  
	\sum_{j=1}^r \left(-z_j,z_j\right)\right)\\
&= \left(0, (\mp x)^{\mu_2}, (\pm x)^{\mu_3},
	\sum_{j=1}^r \left(-z_j,z_j\right)\right)\\
&= f_{r+1}\left((\mp x)^{\mu_2},(\pm x)^{\mu_3},\sum_{j=1}^r z_j\right)\\	
&=f_{r+1}(d_{r+2}(x)), 
\end{align*}
where $z_j$ denotes the image of $x$ under $d_{r+2}$ in $\calM_{3,i_j}$.

The remaining cases to be checked are similar to each other.  Consider the restriction of $f$ to $\calM_{X_1}\oplus\calM_{X_2}$.  Then
$$\xymatrix{
\calM_{X_1}\oplus \calM_{X_2}
	\ar[r]^-{d_{r+2}}\ar[d]_{f_{r+2}=\mathrm{id}} &
\calM_{X_{12}}\oplus\calM_{X_{13}}\oplus\calM_{X_{23}} \oplus \calN
	\ar[d]^{f_{r+1}=\mathrm{id}} \\  
\calM_{X_1}\oplus\calM_{X_2} 
	\ar[r]^-{d_{r+2}}&
\calM_{X_{12}}\oplus\calM_{X_{13}}\oplus\calM_{X_{23}} \oplus \calN\\
}$$
where $\calN$ is a direct sum of modules such that each summand $\calM_{F'}$ is indexed by a state $F'$ that contains $e_3$ and one of $e_1$ or $e_2$.  So $f|_{\calN}=\mathrm{id}$, and thus $f_{r+1}d_{r+2}=d_{r+2}f_{r+2}$.

Next, consider the restriction of $f$ to $\calM_{X_{12}}\oplus\calM_{X_{13}}\oplus\calM_{X_{23}}$. Then
$$\xymatrix{
\calM_{X_{12}}\oplus\calM_{X_{13}}\oplus\calM_{X_{23}} 
	\ar[d]_{f_{r+1}=\mathrm{id}} \ar[r]^-{d_{r+1}}&
	\calM_{X_{123}}\oplus \calN \ar[d]^{f_r=\mathrm{id}} \\
\calM_{X_{12}}\oplus\calM_{X_{13}}\oplus\calM_{X_{23}} 
	\ar[r]^-{d_{r+1}} & \calM_{X_{123}}\oplus\calN \\
}$$
where $\calN$ is a direct sum of modules such that each summand $\calM_{F'}$ is indexed by a state $F'$ that contains at most one of the edges $e_1,e_2$ or $e_3$.  So $f|_{\calN}=\mathrm{id}$, and thus 
$f_{r}d_{r+1} = d_{r+1}f_{r+1}$.

Lastly, consider the restriction of $f$ to $\calM_{X_{123}}$.  Then
$$\xymatrix{
\calM_{X_{123}}\ar[d]_{f_r=\mathrm{id}} \ar[r]^-{d_r}
	& \calN \ar[d]^-{f_{r-1}=\mathrm{id}}\\
\calM_{X_{123}} \ar[r]^-{d_r} 
	& \calN \\
}$$
where $\calN$ is a direct sum of modules such that each summand $\calM_{F'}$ is indexed by a state $F'$ that does not contain any of the edges $e_1,e_2$ or $e_3$.  So $f|_{\calN} = \mathrm{id}$, and thus
$f_{r-1}d_r = d_r f_r$.

Therefore, $f_*$ is a chain map.

It is straightforward to show that $f\circ \iota = \mathrm{id}_{C_*(A\cup B)}$.  It suffices to check this when restricted to each of the $2^{m-3}$ sublattices indexed by a state $X\supseteq T$.  Note that 
$f_{r+2}(\iota(x,y)) = f_{r+2}(x,y,0) = (x,y)$.
$$\xymatrix{
0\ar[r]^-{d_{r+3}} \ar[d]_-\iota& 
	\calM_{X_1}\oplus\calM_{X_2} \ar[r]^-{d_{r+2}} \ar[d]_-\iota& 
	\calM_{X_{12}}\oplus\calM_{X_{13}}\oplus\calM_{X_{23}}
		\ar[r]^-{d_{r+1}} \ar[d]_-\iota& 
	\calM_{X_{123}} \ar[d]_-\iota\\
\calM_X\ar[d]_{f_{r+3}=0} \ar[r]^-{d_{r+3}}& 
	\calM_{X_1}\oplus\calM_{X_2} \oplus \calM_{X_3}	 
		\ar[r]^-{d_{r+2}} \ar[d]_{f_{r+2}} & 
	\calM_{X_{12}}\oplus\calM_{X_{13}}\oplus\calM_{X_{23}}
		\ar[r]^-{d_{r+1}} \ar[d]_{f_{r+1}=\mathrm{id}} & 
	\calM_{X_{123}} \ar[d]_{f_{r}=\mathrm{id}} \\
0\ar[r]^-{d_{r+3}}& 
	\calM_{X_1}\oplus\calM_{X_2} \ar[r]^-{d_{r+2}} & 
	\calM_{X_{12}}\oplus\calM_{X_{13}}\oplus\calM_{X_{23}}
		\ar[r]^-{d_{r+1}} & 
	\calM_{X_{123}} \\
}$$

Finally, we show that $\iota\circ f$ is chain homotopic to the identity.  Again, partition the states of $G$ into $2^{m-3}$ sublattices, each indexed by an $X \supseteq T$, and define $s_*$ on the eight corresponding modules.

Suppose $X=\{e_1,e_2,e_3,e_{i_1},\ldots, e_{i_r}\}$, and without loss of generality, we may assume that $X$ is connected, for otherwise, define $s$ on the component of $X$ containing $T$, then extend $s$ to be identity on the remaining components, and finally induce $s\otimes \mathrm{id}$ to a map on $\fS_n$-modules. 
\begin{enumerate}
\item[(i)] Since $T$ is a triangle, then deleting the edge $e_3$ will not disconnect $X$, and therefore, $\calM_X = \calM_{X_3}$.
If $F=X_3$, let $s_{r+2}(x) = \mathrm{sgn}(d_\varepsilon)x$ for any $x\in X_3$,
where $\mathrm{sgn}(d_\varepsilon)$ is the sign of the per-edge map $d_\varepsilon:X \rightarrow X_3$ in $Q(G)$.
\item[(ii)] If $F$ is any of the seven remaining states, let $s|_F=0$.
\end{enumerate}
The following diagram illustrates $s_*$ when restricted to a sublattice indexed by some $X\supseteq T$.  To compute $ds+sd$, one needs to direct sum all $2^{m-3}$ such sublattices.
$$\xymatrix@!C=110pt{
\calM_X \ar[d]|{f_{r+3}=0} 
		\ar[r]^-{d_{r+3}} & 
	\calM_{X_1}\!\oplus\!\calM_{X_2}\!\oplus\!\calM_{X_3}	 
		\ar[r]^-{d_{r+2}} 
		\ar[d]|{f_{r+2}} 
		\ar[ldd]_<<<<<<<<<<{s_{r+2}} & 
	\calM_{X_{12}}\!\oplus\!\calM_{X_{13}}\!\oplus\!\calM_{X_{23}}
		\ar[r]^-{d_{r+1}} 
		\ar[d]|{f_{r+1}=\mathrm{id}} 
		\ar[ldd]_<<<<<<<<<<{s_{r+1}=0}& 
	\calM_{X_{123}} 
		\ar[d]|{f_{r}=\mathrm{id}} 
		\ar[ldd]_<<<<<<<<<<{s_{r}=0}\\
0\ar[d]_-\iota& 
	\calM_{X_1}\!\oplus\!\calM_{X_2} \ar[d]_-\iota & 
	\calM_{X_{12}}\!\oplus\!\calM_{X_{13}}\!\oplus\!\calM_{X_{23}}
		\ar[d]_-\iota & 
	\calM_{X_{123}} \ar[d]_-\iota \\
\calM_X\ar[r]^-{d_{r+3}}& 
	\calM_{X_1}\!\oplus\!\calM_{X_2}\!\oplus\!\calM_{X_3}	 
		\ar[r]^-{d_{r+2}} & 
	\calM_{X_{12}}\!\oplus\!\calM_{X_{13}}\!\oplus\!\calM_{X_{23}}
		\ar[r]^-{d_{r+1}} & 
	\calM_{X_{123}} \\
}$$

Since the prism operators are almost all zero, it is only necessary to check that $\mathrm{id}-\iota\circ f = ds+sd$ when restricted to $\calM_X$ and $\calM_{X_3}$; the other cases are clear.

Consider the case $F=X$.  
$$\xymatrix@!C=100pt{
\ &
\calM_X \ar[d]^{f_{r+3}=0} 
		\ar[r]^-{d_{r+3}} 
		\ar[ldd]_{s_{r+3}=0} & 
	\calM_{X_1}\oplus\calM_{X_2}\oplus\calM_{X_3}\oplus\calN	 
		\ar[ldd]^{s_{r+2}}\\
 &
0\ar[d]^-\iota
& 
\\
\calO\ar[r]^-{d_{r+4}} &
\calM_X 
	& 
	\\
}$$
where $\calN=\bigoplus_{j=1}^r \calM_{X_{i_j}}$ such that $X_{i_j}\supseteq T$, so $s|_{\calN}=0$ (and by definition $s_{r+3}$ is nonzero only on the summand $\calM_{X_3}$).

Thus, for all $x\in \calM_X$,
$(\mathrm{id}-\iota\circ f)(x) = x$, and
$$(d_{r+4}s_{r+3}+ s_{r+2}d_{r+3})(x) = 0 +  s_{r+2}(\pm x, \mp x, \pm x, n) = \pm (\pm x) = x. $$

As for the case $F=X_3$, 
$$\xymatrix@!C=100pt{
&\calM_{X_3}	 
	\ar[r]^-{d_{r+2}} 
	\ar[d]^-{f_{r+2}} 
	\ar[ldd]_{s_{r+2}}
	& \calM_{X_{13}}\oplus\calM_{X_{23}}\oplus\calN
		\ar[ldd]^{s_{r+1}}\\
&\calM_{X_1}\oplus\calM_{X_2} \ar[d]^-\iota \\
\calM_X \ar[r]^-{d_{r+3}}
	& \calM_{X_1}\oplus\calM_{X_2}\oplus \calM_{X_3}\oplus\calP
}$$
where $\calN=\bigoplus_{j=1}^r \calM_{X_{3,i_j}}$ is a direct sum of modules such that each summand is indexed by a state $F'$ that contains $e_1$ and $e_2$ but not $e_3$, and $\calP=\bigoplus_{j=1}^r \calM_{X_{i_j}}$.  

Note, in particular, that 
$$\xymatrix@R=1pt@C=80pt{
\calM_{X_3}\ar[r]^{s_{r+2}=\mathrm{sgn}(d_\varepsilon)\mathrm{id}} 
	& \calM_X \ar[r]^{\mathrm{sgn}(d_\varepsilon)\mathrm{id}} 
	& \calM_{X_3}\\
x \ar@{|->}[r]
	& \mathrm{sgn}(d_\varepsilon)x \ar@{|->}[r] 
	& \mathrm{sgn}(d_\varepsilon)x = x,
}$$
so it follows that
$$d_{r+3}(s_{r+2}(x))
=d_{r+3}(\mathrm{sgn}(d_\varepsilon)x)
=\left(x, -x, x, \sum_{j=1}^r z_j\right),
$$
where $\varepsilon$ is the edge from $X$ to $X_3$ in $Q(G)$, and $z_j$ is image of $x$ in $\calM_{X_{3,i_j}}$ under $d$.

On the other hand, note that $\calM_{X_{3,i_j}}= \calM_{X_{i_j}}$ because they have the same connected graph components, so the per-edge maps 
$\calM_{X_3}\rightarrow \calN$ and $\calM_X \rightarrow \calP$ are identical.  As well, it follows from this that $s$ is the identity map when restricted to $\calN$.
Moreover, if $\eta_j$ is the edge from $X_{i_j}$ to $X_{3,i_j}$ in $Q(G)$, then $\mathrm{sgn}(d_{\eta_j}) = - \mathrm{sgn}(d_\varepsilon)$, since $i_j>3$ by assumption.  Also, we mention that by definition,
$s$ is zero on $\calM_{X_{13}}\oplus \calM_{X_{23}}$, so
$$s_{r+1}(d_{r+2}(x)) 
=s_{r+1}\left((-x)^{\mu_2}, x^{\mu_3}, \sum_{j=1}^r z_j\right)
=\left(0,0,0, -\sum_{j=1}^r z_j\right),
$$
and therefore,
$$\left(d_{r+3}s_{r+2}+s_{r+1}d_{r+2}\right)(x) = (x,-x,x,0).$$

Lastly, for all $x\in\calM_{X_3}$, $(\mathrm{id}-\iota\circ f)(x) = (0,0,x,0) - (-x,x,0) = (x,-x,x,0)$.  So $\mathrm{id}-\iota\circ f = ds +sd$ when restricted to $\calM_{X_3}$.

This completes the proof that $\mathrm{id}_{C_*(G)}$ is homotopic to $\iota\circ f$.  Therefore, $H_*(G) \cong H_*(A\cup B)$. \qed

\begin{lemma} \label{lem.ses}
Let $G$ be a graph where the edges $e_1, e_2, e_3\in E(G)$ form a triangle.  Let $A = G-e_1$ and $B=G-e_2$.
For $i\geq0$, there are short exact sequences
$$\xymatrix{ 0 \ar[r]& C_i(A\cap B) \ar[r]^-\varphi& C_i(A) \oplus C_i(B) \ar[r]^-\psi& C_i(A\cup B) \ar[r]&0
}$$
where $\varphi(x) = (x, -x)$, and $\psi(a,b)= a+b$.
\end{lemma}
\proof

We have $\ker\varphi = 0$, since $\varphi(x)=(x,-x) =(0,0)$ implies $x=0$.  Next, since $(\psi\circ\varphi)(x) = \psi(x,-x) = x-x=0$, then $\im \varphi \subseteq \ker \psi$. Also, $\ker \psi \subseteq \im \varphi$, since if 
$\psi(a,b) = a+b=0$, then $a=-b\in C_i(B)$, so $a \in C_i(A\cap B)$, and $(a,b) = \varphi(a) \in \im \varphi$.
Lastly, $\im\psi = C_i(A\cup B)$ since $\psi$ acts as identity when restricted to either $C_i(A)$ or $C_i(B)$.
\qed

\begin{theorem}\label{thm.les} Let $G$ be a graph where the edges $e_1, e_2, e_3\in E(G)$ form a triangle.  Let $A = G-e_1$ and $B=G-e_2$.
Then there is a long exact sequence 
$$\xymatrix{\cdots\ar[r]&
H_i(A \cap B) \ar[r]^-{\varphi} &
H_i(A) \oplus H_i(B) \ar[r]^-{\psi} &
H_i(G) \ar[r]^-{\delta} & \cdots \ar[r]&
H_0(G) \ar[r] & 0.
} $$
\end{theorem}
\proof  The short exact sequences of chain complexes in Lemma~\ref{lem.ses} induce a long exact sequence in homology.  
By Lemma~\ref{lem.modchain}, since $H_*(A\cup B)\cong H_*(G)$, then the result follows.\qed

\begin{corollary}
$$\Frob_G(q,1) = \Frob_{G-e_1}(q,1) + \Frob_{G-e_2}(q,1) -\Frob_{G-e_1-e_2}(q,1).$$
\end{corollary}
\proof
The maps in Theorem~\ref{thm.les} are degree-preserving, so applying the Frobenius characteristic map gives the result.\qed

\begin{remark} Further setting $q=1$ in the above Corollary recovers the result in Theorem~\ref{thm.dd}.
\end{remark}

\begin{example}
Consider $G=K_3$.
Theorem~\ref{thm.dd} gives
$$
\begin{array}{ccccccc}
\beginpicture
	\setcoordinatesystem units <.7cm,.7cm>         
	\setplotarea x from 0 to 1, y from 0 to 2    
	\multiput{$\bullet$} at 0 0 1 0 0.5 0.867 /
	\put{\tiny$e_1$}[c] at 1.1 .5 \put{\tiny$e_2$}[c] at -.1 .5 \put{\tiny$e_3$} at .5 -.3
	\put{$G$} at .5 -1
	\plot 0 0 1 0 0.5 .867 0 0 /
	\endpicture
	&=
	&\beginpicture
	\setcoordinatesystem units <.7cm,.7cm>         
	\setplotarea x from 0 to 1, y from 0 to 2    
	\multiput{$\bullet$} at 0 0 1 0 0.5 0.867 /
	\put{\tiny$e_2$}[c] at -.1 .5 \put{\tiny$e_3$} at .5 -.3
	\put{$A$} at .5 -1
	\plot 0.5 .867 0 0 1 0 /
	\endpicture
	&+
	&\beginpicture
	\setcoordinatesystem units <.7cm,.7cm>         
	\setplotarea x from 0 to 1, y from 0 to 2    
	\multiput{$\bullet$} at 0 0 1 0 0.5 0.867 /
	\put{\tiny$e_1$}[c] at 1.1 .5 \put{\tiny$e_3$} at .5 -.3
	\put{$B$} at .5 -1
	\plot 0 0 1 0 0.5 .867 /
	\endpicture
	&-
	&\beginpicture
	\setcoordinatesystem units <.7cm,.7cm>         
	\setplotarea x from 0 to 1, y from 0 to 2    
	\multiput{$\bullet$} at 0 0 1 0 0.5 0.867 /
	\put{\tiny$e_3$} at .5 -.3
	\put{$A\cap B$} at .5 -1
	\plot 0 0 1 0 /
	\endpicture
\\ \\
6s_{111}
	&=
	&(s_{21}+4s_{111})
	&+
	&(s_{21}+4s_{111})
	&-
	&(2s_{21}+2s_{111})
\end{array}
$$

To condense notation, let $H_*(Z)=H_*(A\cap B)$ and $H_*(Y) = H_*(A)\oplus H_*(B)$.  Using the computations from Sections~\ref{eg.kthree}, \ref{eg.smalltrees} and Section~\ref{eg.ktwo} combined with Proposition~\ref{prop.disjointunion}, the long exact sequence for $G=K_3$ is
$$
\xymatrix@C-1.35pc@R-1pc 
{H_2Z \ar[r]& H_2Y\ar[r]& H_2G\ar[r]&
	H_1Z \ar[r]& H_1Y\ar[r]& H_1G\ar[r]&
	H_0Z \ar[r]& H_0Y\ar[r]& H_0G\ar[r]& 0\\
0\ar[r]
	&\hbox{\scriptsize$2\mathcal{S}^{111}$}\ar[r]
	&\hbox{\scriptsize$2\mathcal{S}^{111}$}\ar[r]
	&0\ar[r]& 0\ar[r]& 0\ar[r]& 
	0\ar[r]& 0\ar[r]& 0\ar[r]&
	0 
\\
0\ar[r]& 0\ar[r]\ar[u]& 0\ar[r]\ar[u]
	&\hbox{\scriptsize$\mathcal{S}^{21}\!\!\!+\!\mathcal{S}^{111}$}
		\ar[r]\ar[u]
	&\hbox{\scriptsize$2\mathcal{S}^{21}\!\!\!+\!4\mathcal{S}^{111}$}
		\ar[r]\ar[u]
	&\hbox{\scriptsize$\mathcal{S}^{21}\!\!\!+\!3\mathcal{S}^{111}$}
		\ar[r]\ar[u]
	&0\ar[r]\ar[u]
	&0\ar[r]\ar[u]
	& 0\ar[r]\ar[u]
	&0 
\\
0\ar[r]& 0\ar[r]\ar[u]& 0\ar[r]\ar[u]&
	0\ar[r]\ar[u]& 
	0\ar[r]\ar[u]
	&\hbox{\scriptsize$\mathcal{S}^{21}$}\ar[r]\ar[u]
	&\hbox{\scriptsize$\mathcal{S}^{21}\!\!\!+\!\mathcal{S}^{111}$}
		\ar[r]\ar[u]
	&\hbox{\scriptsize$2\mathcal{S}^{111}$}\ar[r]\ar[u]
	&\hbox{\scriptsize$\mathcal{S}^{111}$}\ar[r]\ar[u]
	&0
}
$$

In terms of the Frobenius series $\Frob_G(q,1)$, 
$$\begin{array}{ccccccc}
\Frob_{K_3}(q,1) &=& \Frob_A(q,1) &+& \Frob_B(q,1) &-& \Frob_{A\cap B}(q,1)\\
\hline
2q^2 s_{111}  
	&=
	&q^2 s_{111} 
	&+
	&q^2 s_{111} 
	&-
	&0 
\\
q(s_{21}+3s_{111})  
	&=
	&q(s_{21}+2s_{111}) 
	&+
	&q(s_{21}+2s_{111})
	&-
	&q(s_{21}+s_{111}) 
\\
s_{111}-s_{21}  
	&=
	&s_{111} 
	&+
	&s_{111} 
	&-
	&(s_{21}+s_{111})
\end{array}
$$

\end{example}

\section{Examples}\label{sec.computations}

We include computations of chromatic symmetric homology for several small graphs.

\subsection{Edge-less graphs}\label{eg.edgeless}
  
Let $G$ be the edge-less graph on $n$ vertices.  In this case, $G$ only has one state, and the $\fS_n$-module associated to the state is 
$$\Ind_{\fS_1\times\cdots\times \fS_1}^{\fS_n} \left( \calS^{1}\otimes \cdots \otimes \calS^{1}\right)=\bigoplus_{\lambda\vdash n} f^\lambda \calS^\lambda,$$
where $f^\lambda=\chi^\lambda(1^n)=\dim(\calS^\lambda)$ is the number of standard Young tableaux of shape $\lambda$.
Thus $H_{00}= \bigoplus_{\lambda\vdash n} f^\lambda \calS^\lambda$, and the bigraded Frobenius series 
$\Frob_G(q,t) = \sum_{\lambda\vdash n} f^\lambda s_\lambda = p_1^n = X_G$
is independent of $q$ and $t$.

\subsection{A single edge}\label{eg.ktwo} 
Let $G=K_2$.  The diagram of states together with the signed per-edge map is 
$$ 
\beginpicture
	\setcoordinatesystem units <4cm,1cm>         
	\setplotarea x from 0 to 1, y from -.75 to .5    

\arrow <8pt> [.2,.67] from 0.2 0 to .8 0    \put{$d_{*}$} at .5 .3

\put{$\beginpicture
	\setcoordinatesystem units <.7cm,.7cm>         
	\setplotarea x from 0 to 1, y from 0 to 0    
	\put{$\bullet$} at 0 0 \put{$\bullet$} at 1 0 
	\plot 0 0 1 0 /
	\put{\tiny1} at .5 -.5
	\endpicture$} at 0 0	

\put{$\beginpicture
	\setcoordinatesystem units <.7cm,.7cm>         
	\setplotarea x from 0 to 1, y from 0 to 0    
	\put{$\bullet$} at 0 0 \put{$\bullet$} at 1 0
	\put{\tiny0} at .5 -.5
	\endpicture$} at 1 0
	
\endpicture
$$
The graded $\fS_2$-modules corresponding to the states are:
$$ 
\xymatrix{
\calS^{11} \\
\calS^{2} \ar[r]^-{d_{10}} \ar@{.>}[u]
	& \Ind_{\fS_1\times \fS_1}^{\fS_2} 
	\left(\calS^{1}\otimes \calS^{1}\right) \\
}
$$
Recall that $\fS_2 = \langle s_1 \mid s_1^2 =e\rangle$.  Then
$\Ind_{\fS_1\times \fS_1}^{\fS_2}\left(\calS^1\otimes \calS^1\right) 
\cong \bbC[\fS_2]
\cong \calS^2 \oplus \calS^{11} 
= \spn_\bbC\{e,s_1\}$. 
Choosing $\calS^2 = \spn_\bbC\{ e+s_1\}$ and $\calS^{11}=\spn_\bbC\{e-s_1\}$, then 
$$d_{10}: \spn_\bbC\{e+s_1\} \rightarrow \spn_\bbC\{e,s_1\}
$$
is simply the inclusion map.  Thus $\ker d_{11}\cong\calS^{11}$, $\im d_{10} \cong \calS^2$, and $\ker d_{00}\cong \calS^2\oplus \calS^{11}$. In summary, the chromatic symmetric homology of graded $\fS_2$-modules for  $G=K_2$ is
$$
\begin{array}{c|cc}
1&
	\calS^{11}&\\
0&	
	0&\calS^{11}\\
\hline
& H_1(K_2)  & H_0(K_2)  
\end{array}
$$
and the bigraded Frobenius series is 
$\Frob_{K_2}(q,t)=qts_{11}+s_{11} = (1+qt)s_{11}.$

\subsection{Some small trees}
While the chromatic polynomial $\chi_G$ does not distinguish trees, it is conjectured that the chromatic symmetric function is a complete invariant for trees.  In this example, we state results in homology for trees up to four vertices.

\begin{example} \label{eg.smalltrees} 
Consider the two trees on four vertices.  Let
$$
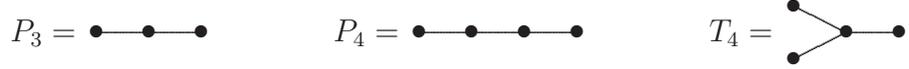
\beginpicture
	\setcoordinatesystem units <.7cm,.7cm>         
	\setplotarea x from 0 to 2, y from 0 to 1    
	\multiput{$\bullet$} at 0 .5 1 .5 2 .5 /
	\plot 0 .5 2 .5 /
	\put{$P_3=$} at -1 .5
\endpicture
\qquad
\qquad
\beginpicture
\setcoordinatesystem units <.7cm,.7cm>         
\setplotarea x from 0 to 2, y from 0 to 1    
\multiput{$\bullet$} at 0 0.5 1 0.5 2 0.5 3 0.5 / 
\plot 0 0.5 3 0.5 /
\put{$P_4=$} at -1 0.5
\endpicture
\qquad
\qquad
\beginpicture
\setcoordinatesystem units <.7cm,.7cm>         
\setplotarea x from 0 to 2, y from 0 to 1    
\multiput{$\bullet$} at 0 0 0 1 1 0.5 2 0.5 /
\plot 0 1 1 0.5 0 0 /
\plot 1 0.5 2 0.5 /
\put{$T_4=$} at -1 0.5
\endpicture
.$$

The chromatic symmetric homology of graded $\fS_n$-modules for each of these graphs is

$$\begin{array}{c|ccc}
2 
	& \calS^{111} & &  \\
1 
	& 0
	& \calS^{21}\oplus\left(\calS^{111}\right)^{\oplus2} & \\
0 
 	& 0 
	& 0
	& \calS^{111} \\
\hline 
& H_2(P_3)& H_1(P_3)& H_0(P_3)
\end{array}$$
\bigskip

$$\begin{array}{c|cccc}
3
	& \calS^{1111} \\
2 
	& 0
	&\calS^{22}\oplus (\calS^{211})^{\oplus2} 
		\oplus (\calS^{1111})^{\oplus3}  \\
1 
	& 0
	& 0
	& \calS^{22}\oplus (\calS^{211})^{\oplus2} 
		\oplus (\calS^{1111})^{\oplus3}\\
0 
 	& 0 
	& 0
	& 0
	& \calS^{1111} \\
\hline 
& H_3(P_4)& H_2(P_4)& H_1(P_4)& H_0(P_4)\\
\end{array}
$$
\bigskip
and
$$\begin{array}{c|cccc}
3
	& \calS^{1111}\\
2 
	& 0
	& (\calS^{211})^{\oplus2} \oplus (\calS^{1111})^{\oplus3}\\
1 
	& 0
	& \calS^{22} 
	& \calS^{31}\oplus\calS^{22}\oplus(\calS^{211})^{\oplus3} \oplus (\calS^{1111})^{\oplus3}\\
0 
	& 0
 	& 0 
	& \calS^{22}
	& \calS^{1111} \\
\hline 
& H_3(T_4) & H_2(T_4)& H_1(T_4)& H_0(T_4)\\
\end{array}
$$
\bigskip

Their Frobenius series are
\begin{align*}
\Frob_{P_3}(q,t)&=qts_{21} +(1+qt)^2s_{111},\\
\Frob_{P_4}(q,t)&=qt(1+qt)s_{22}+2qt(1+qt)s_{211}+(1+qt)^3s_{1111},\\ 
\Frob_{T_4}(q,t)&=qts_{31}-t(1-q+qt)s_{22} +qt(3+2qt)s_{211}+(1+qt)^3s_{1111}. 
\end{align*}

\end{example}


\end{document}